\documentclass[a4paper, 12pt, reqno]{amsart}
\usepackage{anyfontsize}
\usepackage{amsmath} 
\usepackage{amssymb}
\usepackage{amsthm}
\usepackage{mathrsfs}
\usepackage{geometry}
\usepackage{tikz}
\usepackage{tikz-cd}
\usepackage{enumitem}
\usepackage{booktabs}
\usepackage{subfigure}
\usepackage[graphicx]{realboxes}
\usepackage{caption}
\usetikzlibrary{commutative-diagrams}
\usepackage{hyperref}
\theoremstyle{definition}
\newtheorem{theorem}{Theorem}[section]
\newtheorem{corollary}[theorem]{Corollary}
\newtheorem{lemma}[theorem]{Lemma}

\newtheorem{proposition}[theorem]{Proposition}

\newtheorem{remark}[theorem]{Remark}
\newtheorem{example}[theorem]{Example}
\newtheorem{construction}[theorem]{Construction}
\newtheorem{condition}[theorem]{Condition}

\allowdisplaybreaks[2]
\numberwithin{equation}{section}

\DeclareMathOperator{\id}{id}
\DeclareMathOperator{\im}{im}
\DeclareMathOperator{\GL}{GL}
\DeclareMathOperator{\Hom}{Hom}
\DeclareMathOperator{\supp}{supp}
\DeclareMathOperator{\Z}{\mathcal{Z}}
\DeclareMathOperator{\ZZ}{\mathbb{Z}_2}

\begin{document}
\title[EQUIVARIANT BORDISM CLASSIFICATION]{EQUIVARIANT BORDISM CLASSIFICATION OF FIVE-DIMENSIONAL $(\ZZ)^3$-MANIFOLDS WITH ISOLATED FIXED POINTS}
\author{Yuanxin Guan}
\address{School of Mathematical Sciences, Fudan University, Shanghai, 200433, China}
\email{yxguan21@m.fudan.edu.cn}
\author{Zhi L\"u}
\address{School of Mathematical Sciences, Fudan University, Shanghai, 200433, China}
\email{zlu@fudan.edu.cn}

\begin{abstract}
	Denote by $\Z_5((\ZZ)^3)$ the group, which is also a vector space over $\ZZ$, generated by equivariant unoriented bordism classes of all five-dimensional closed smooth manifolds with effective smooth $(\ZZ)^3$-actions fixing isolated points. We show that $\dim_{\ZZ} \Z_5((\ZZ)^3) = 77$ and determine a basis of $\Z_5((\ZZ)^3)$, each of which is explicitly chosen as the projectivization of a real vector bundle. Thus this gives a complete classification up to equivariant unoriented bordism of all five-dimensional closed smooth manifolds with effective smooth $(\ZZ)^3$-actions with isolated fixed points.
\end{abstract}




\maketitle

\section{Introduction}
Let $\Z_n((\ZZ)^k)$ be the group generated by the equivariant unoriented bordism classes of $n$-dimensional closed smooth manifolds with effective $(\ZZ)^k$ actions fixing isolated points, and $\Z_*((\ZZ)^k) = \oplus_{n\geqslant 0} \Z_n((\ZZ)^k)$, which is a graded algebra over $\ZZ$. 
In \cite{ConnerFloyd}, Conner and Floyd initialized the study of $\Z_*((\ZZ)^k)$, and introduced the Conner-Floyd representation algebra $R_*((\ZZ)^k)$, isomorphic to the polynomial ring over $\ZZ$ with $2^k$ indeterminates (see also \cite{Conner1979, Stong} and Section \ref{section_ConnerFloyd_repre} for definition). They carried out an elaborate analysis of actions of $\ZZ\times \ZZ$ with isolated fixed points via the graded algebra homomorphism $\phi_*: \Z_*((\ZZ)^k)\to R_*((\ZZ)^k)$, by assigning the tangent representations at all fixed points to every equivariant bordism class. 
The classical theorem of Stong in \cite{Stong} implied that $\phi_*$ is a monomorphism. Thus
$\Z_*((\ZZ)^k)$ is isomorphic to $\im\phi_*$, which prompts us to study the representation theory of fixed points. This provides an approach to the study of $\Z_*((\ZZ)^k)$.
For $k=1,2$, we have known from Conner-Floyd's work \cite[Theorems 25.1, 31.1 and 31.2]{ConnerFloyd} that $\Z_*((\ZZ))$ is trivial, and $\Z_*(\ZZ\times \ZZ)$ is the polynomial algebra over $\ZZ$ with one variable which is the bordism class of the real projective plane with standard $\ZZ\times \ZZ$-action. To the best of our knowledge, the ring structure of $\Z_*((\ZZ)^k)$ for larger $k$ is not determined yet, in view of the exponential growth of the number of generators in $R_*((\ZZ)^k)$. Over the past 15 years, the primary focus was on the following questions:
\begin{enumerate}
       \item[({\bf Q1})]
       {\em What polynomials in $\mathcal{R}_*((\ZZ)^k)$ arise as fixed point
data of $(\ZZ)^k$-manifolds?}
		\item[({\bf Q2})] {\em Is there a preferred representative in each equivariant bordism class of $\Z_n((\ZZ)^k)$?}
     \item[({\bf Q3})]{\em What is the dimension of $\Z_n((\ZZ)^k)$ as a vector space over $\ZZ$ for every $n$ and $k$?} Note that $\Z_n((\ZZ)^k)$ is trivial for $0<n<k$ by the effectiveness of $(\ZZ)^k$-actions. 
\end{enumerate}

In \cite{ChenLuTan, Lu2009, LuTan2014}, Chen, L\"u and Tan completely answered the above three questions for the case $n=k$. Specifically, they obtained that 
\begin{itemize}
 \item The image $\phi_*( \mathcal{Z}_n((\ZZ)^n))$ is characterized by the differential on the dual algebra of the Conner--Floyd representation algebra.
 \item Each class in $\mathcal{Z}_n((\ZZ)^n)$ contains a small cover as its representative, and in particular, $\mathcal{Z}_n((\ZZ)^n)$
 is generated by the classes of generalized real Bott $n$-manifolds (i.e., small covers over products of simplices).
 Here, small covers, 
served as the topological counterparts of real toric varieties, were introduced and studied by Davis and Januszkiewicz in their seminar work~\cite{smallcover}. 
 \item $\dim_{\ZZ}\Z_n((\ZZ)^n)=(-1)^n + \sum_{i=0}^{n-1} (-1)^{n-1-i} \frac{(2^n-2^i)\cdots (2^n-2^0)}{(i+1)!}$.
\end{itemize}
As noted in \cite{ChenLuTan}, the research work on $\mathcal{Z}_n((\ZZ)^n)$ demonstrated profound connections among various theories, including equivariant bordism, ordinary bordism of spaces, homology, representation theory, Davis-Januszkiewicz theory, and matroid theory. Notably, within the framework of equivariant bordism, the Davis-Januszkiewicz theory proves applicable to a broader class of objects beyond just small covers. 

The GKM theory, established by Goresky, Kottwitz and MacPherson \cite{GKMtheory} for algebraic varieties with complex torus actions, demonstrates that the information of tangent representations at fixed points of actions can be recoded on labeled graphs.
Recently, Li, L\"u and Shen in \cite{LiLuShen_2025arxiv} brought this characteristic of the GKM theory into the study of $\Z_*((\ZZ)^k)$, forming a relationship between equivariant bordism and the GKM theory. They gave two characterizations of all polynomials of the subalgebra $\phi_*(\Z_*((\ZZ)^k))$ in terms of $(\ZZ)^k$-labeled graphs and $(\ZZ)^k$-representations, respectively, answering the question $({\bf Q1})$ from the perspectives of combinatorics and $(\ZZ)^k$-representation theory. This makes it feasible to determine whether a given polynomial of $R_*((\ZZ)^k)$ is the image of some $(\ZZ)^k$-manifold in $\Z_*((\ZZ)^k)$ without finding the preimage.
As an application, it was shown by a nontrivial argument in \cite{LiLuShen_2025arxiv} that $\dim_{\ZZ} \Z_4((\ZZ)^3) = 32$ with explicit geometric generators constructed, yielding an exhaustive categorization up to equivariant unoriented bordism of all 4-dimensional closed smooth manifolds with effective $(\ZZ)^3$-actions fixing a finite set.
This breaks the case of $n=k$, obtaining the first result on the study of $\Z_n((\ZZ)^k)$ with $n>k\geqslant 3$. Meanwhile, the characterization result of Li-L\"u-Shen in \cite{LiLuShen_2025arxiv} appears potentially applicable to uncovering the structure of $\Z_n((\ZZ)^k)$.

In this paper, we focus on the case of $n=5$ and $k=3$, and successfully settle this case. 
We shall see that every nonzero polynomial in $\phi_*(\Z_5((\ZZ)^3))$ has a monomial in its support divided by the square of some variable, which leads us to the determination of a quantity of polynomials in $\phi_*(\Z_5((\ZZ)^3))$, thus obtaining the explicit solution of the question (\textbf{Q1}) (see Theorem~\ref{generating_set}). Moreover, we derive the dimension of $\Z_5((\ZZ)^3)$ that resolves the question (\textbf{Q3}) in this special case.
\begin{theorem}[Theorem \ref{thm_dimension}]
$\dim_{\ZZ}\Z_5((\ZZ)^3)=77$.
\end{theorem}

Next, we address the question ({\bf Q2}) in the case $n=5$ and $k=3$. 
Equivalently, we shall just construct the required geometric generators of $\Z_5((\ZZ)^3)$. 
Since suitable real Milnor hypersurfaces $H_{m, n}$ can be used as generators of the unoriented bordism ring $\Omega_*^O$, we examine whether the Milnor hypersurfaces with certain $(\ZZ)^k$-actions can be our desired representatives. Based on the existing results in \cite{BasuMukSar2014}, a partial solution of the question ({\bf Q2}) is settled. 
After that, building on the work of \cite{Bukhshtaber1998, Lu2009, LuTan2014}, we start our consideration with $n$-dimensional small covers, and attempt to pull back the actions of $(\ZZ)^n$ to actions of $(\ZZ)^k$ through the natural monomorphism $(\ZZ)^k\to (\ZZ)^n$ for $k<n$. A complete answer to the question ({\bf Q2}) is within reach. For the statement of our result, we first introduce some notation.
A manifold $M$ with a $G$-action is often written as $(\varphi, M)$, where $\varphi$ is the map $G\times M\to M$ that defines the action of $G$ on $M$, or just $(G, M)$ if the map that defines the action is understood.
Given a $G$-manifold $(\varphi, M)$ and an automorphism $\sigma$ of $G$, denote by $\varphi_\sigma$ the new action $\varphi\circ (\sigma\times \id_M)$ of $G$ on $M$. 

\begin{theorem}[Theorem \ref{thm_representative}] \label{main_thm2}
	There exist two actions $\varphi_1$ and $\varphi_2$ of $(\ZZ)^3$ on the real Milnor hypersurface $H_{2, 4}$, and actions $\varphi_3$ and $\varphi_4$ of $(\ZZ)^3$ on two 2-staged generalized real Bott manifolds $M_1$ and $M_2$, respectively, such that the set
	\begin{align*}
		\{(\varphi_{1, \sigma}, H_{2, 4}), (\varphi_{2, \sigma}, H_{2, 4}), (\varphi_{3, \sigma}, M_1), (\varphi_{4, \sigma}, M_2): \sigma \text{ is an automorphism of } (\ZZ)^3\}
	\end{align*}
	spans the vector space $\Z_5((\ZZ)^3)$ over $\ZZ$.
\end{theorem}

The actions $\varphi_i$, $i=1, 2, 3, 4$, are constructed in \ref{construct_Milnor_hypersurface}, \ref{construct_small_cover1} and \ref{construct_small_cover2}, respectively.
In particular, since the Milnor hypersurfaces and 2-staged generalized real Bott manifolds are all projectivizations of real vector bundles over real projective spaces (\cite{ChoiMasudaSuh2010, KurokiLu2016, Milnor1965}), Theorem \ref{main_thm2} implies that every such five-dimensional $(\ZZ)^3$-manifold is equivariantly bordant to the projectivization of a real vector bundle over the disjoint union of real projective spaces with certain action of $(\ZZ)^3$.

The paper is organized as follows. 
In Section \ref{section_ConnerFloyd_repre}, we review some work of Conner-Floyd and Stong, which shows that the formal sum of all tangent representations at fixed points is a complete invariant of $(\ZZ)^k$-equivariant bordism.
In Section \ref{section_graph_of_action}, we overview the construction of a $(\ZZ)^k$-labeled graph from an effective $(\ZZ)^k$-manifold with finite fixed points and some properties of the graph, and state the characterization of polynomials in the image of $\phi_*$.
Several basic observations of $\phi_*(\Z_n((\ZZ)^k)$ are delivered in Section \ref{section_dimension}, and the explicit solutions of the questions ({\bf Q1}) and ({\bf Q3}) for $\phi_*(\Z_5((\ZZ)^3)$ are totally determined.
In Section \ref{section_representative}, we settle the question ({\bf Q2}) for the case in which $n=5$ and $k=3$.
Throughout this paper, all manifolds are closed and smooth, and all group actions on manifolds are smooth if no otherwise is stated. 

\section{Conner-Floyd representation algebra and Stong homomorphism} \label{section_ConnerFloyd_repre}
In \cite{ConnerFloyd}, to consider actions of $\ZZ\times \ZZ$ with isolated fixed points, Conner and Floyd defined a representation algebra of all finite dimensional real representations of this group. Let $G$ be a compact Lie group. Denote by $R_n(G)$ the vector space over $\ZZ$ freely generated by all isomorphic classes of $n$-dimensional real representations of $G$, and $R_*(G) = \oplus_{n\geqslant 0} R_n(G)$. Then $R_*(G)$ is equipped with a structure of commutative graded algebra over $\ZZ$ with unit under the direct sum of representations, called the \textit{Conner-Floyd representation algebra}.
Thus $R_*(G)$ is a polynomial algebra over $\ZZ$ whose variables correspond one-to-one to the irreducible representation classes.
In particular, for $G = (\ZZ)^k$, since there is a natural one-to-one correspondence between the set of irreducible representation classes of $(\ZZ)^k$ and the set $\Hom_{\ZZ}((\ZZ)^k, \ZZ)$, the algebra $R_*((\ZZ)^k)$ is isomorphic to the polynomial algebra over $\ZZ$ generated by elements in $\Hom_{\ZZ}((\ZZ)^k, \ZZ)$. From now on, by abuse of notation, we shall denote the formal addition in $R_*((\ZZ)^k)$ by the same symbol as the addition in $\Hom_{\ZZ}((\ZZ)^k, \ZZ)$. 

The significance of the Conner-Floyd representation algebra lies in its ability to describe the information on a fixed point of a group action.
Suppose that a group $G$ acts on an $n$-dimensional manifold $M$ and that $x\in M$ is a fixed point of this action. Then there is an $n$-dimensional real representation of $G$ on the tangent space $T_xM$ at $x$, i.e. a group homomorphism $G\to \GL(T_xM)$, where $\GL(V)$ is the general linear group of a real vector space $V$. We denote the isomorphism class of this representation by $X(x)$, called the \textit{tangent representation at $x$}. When $G$ is the group in which we are interested, the tangent representation at a fixed point of an action was thoroughly studied. Before giving a statement, it is convenient to fix some conventions and notation.
Every element $g$ in $(\ZZ)^k$ in our paper is written as a column vector $g = (g_1, \cdots, g_k)^\mathsf{T}$, where $g_i\in \ZZ=\{0, 1\}$, and by a natural isomorphism $\Hom_{\ZZ}((\ZZ)^k, \ZZ)\cong (\ZZ)^k$, a linear function $(\ZZ)^k\to \ZZ$ can be also regarded as a column vector in the context. Denote by $\rho_0$ the trivial function in $\Hom_{\ZZ}((\ZZ)^k, \ZZ)$, and $\rho_i$ the function $(0, \cdots, 1, \cdots, 0)^\mathsf{T}$ in $\Hom_{\ZZ}((\ZZ)^k, \ZZ)$ where 1 is in the $i$-th position, for $1\leqslant i\leqslant k$. 

\begin{proposition}[\cite{ConnerFloyd}] \label{tangent_representation}
	Suppose that $(\ZZ)^k$ acts on an $n$-dimensional manifold $M$ and $x$ is a fixed point. Then in $R_*((\ZZ)^k)$, 
	\begin{align*}
		X(x) = \prod_{\rho\in \Hom_{\ZZ}((\ZZ)^k, \ZZ)} \rho^{d(\rho)},
	\end{align*}
	such that $\sum d(\rho)=n$, $d(\rho_0)$ is the dimension of the component of the fixed point set $M^{(\ZZ)^k}$ containing $x$, and $d(\rho_0)+d(\rho)$ is the dimension of the component of $M^{\ker \rho}$ containing $x$ for all $\rho\ne\rho_0$.
\end{proposition}

\begin{remark}
	The original form of this result in \cite{ConnerFloyd} was proved for the case where $k=2$ and the fixed points are isolated, and it still holds without these restrictive conditions. 
\end{remark}

This proposition indicates how to determine the tangent representation at a fixed point under the subgroup action. Notice that a subgroup of $(\ZZ)^k$ is also a subspace of $(\ZZ)^k$ as a vector space.

\begin{corollary} \label{subgroup_tangent_repre}
	With the notation of Proposition \ref{tangent_representation}, if $H$ is a subgroup of $(\ZZ)^k$, then there is an action of $H$ on $M$ induced from the action of $(\ZZ)^k$, and the tangent representation at $x$ in the subgroup action, denoted by $X_H(x)$, is 
	\begin{align*}
		X_H(x) = \prod_{\rho\in \Hom_{\ZZ}((\ZZ)^k, \ZZ)} \left(\rho + H^\perp\right)^{d(\rho)},
	\end{align*}
	where $H^\perp = \{(g_1, \cdots, g_k)^\mathsf{T}\in (\ZZ)^k: \sum g_i h_i = 0 \text{ for all } (h_1, \cdots, h_k)^\mathsf{T}\in H\}$ is the orthogonal complement of $H$ in $(\ZZ)^k$.
\end{corollary}

Note that by the identification $\Hom_{\ZZ}((\ZZ)^k, \ZZ)\cong (\ZZ)^k$, $\rho+H^\perp\in (\ZZ)^k/H^\perp\cong H\cong \Hom_{\ZZ}(H, \ZZ)$, and so $\rho+H^\perp$ is an irreducible representation of $H$. 

A group action is often significantly shaped by its fixed point set.
For example, for an involution of a manifold $M$, if the fixed point set is a connected submanifold whose Stiefel-Whitney classes all vanish, then $M$ bounds \cite{ConnerFloyd}, and for a $T^k$- or $(\mathbb{Z}_p)^k$-space $X$ with $p$ prime under certain conditions, the localized homomorphism between the equivariant cohomology rings induced by the inclusion from the fixed point set to $X$ is an isomorphism \cite{HsiangWuyi}.
In our setting, the phenomenon still occurs. There is a natural graded algebra homomorphism $\phi_*: \Z_*((\ZZ)^k)\to R_*((\ZZ)^k)$ defined by
\begin{align*}
	\phi_*([(\ZZ)^k, M]) = \sum_{x\in M^{(\ZZ)^k}} X(x).
\end{align*}

\begin{theorem}[\cite{Stong}] \label{Stong_inj}
	The homomorphism $\phi_*$ is a monomorphism.
\end{theorem}

The theorem indicates that the equivariant bordism class of a $(\ZZ)^k$-manifold with finite fixed points is decided by the tangent representations at fixed points, and the problem of determining the algebraic structure of $\Z_*((\ZZ)^k)$ is converted to determine the algebraic structure of the subalgebra $\im \phi_*$.

\section{Graphs of actions and characterization of the image of \texorpdfstring{$\phi_*$}{}} \label{section_graph_of_action}
We have seen that the formal sum of all tangent representations at fixed points is a complete invariant of equivariant bordism classes of closed smooth $(\ZZ)^k$-manifolds having finite fixed points, and if such a manifold has no fixed point, it must bounds equivariantly.
In this section, we review the construction of a $(\ZZ)^k$-labeled graph from an effective $(\ZZ)^k$-manifold fixing nonempty finite points, such that the graph records all information of fixed points. See \cite{LiLuShen_2025arxiv, Lu2008} for more details.
First we need to clarify that a graph in the paper is an undirected multigraph without loops.

Let $((\ZZ)^k, M^n)$ be a manifold with an effective $(\ZZ)^k$-action fixing nonempty isolated points, and $\rho$ be an arbitrary nontrivial irreducible real representation of $(\ZZ)^k$. Then $\ker\rho$ is a subgroup of $(\ZZ)^k$ of rank $k-1$. Let $C_\rho$ be a connected component of the fixed point set of the action of $\ker\rho$ such that $C_\rho$ is of positive dimension and contains a fixed point in $M^{(\ZZ)^k}$. Moreover, there is an action of $(\ZZ)^k/\ker\rho$ on $C_\rho$ defined by $(g+\ker\rho)\cdot y = gy$. Thus $C_\rho$ is a $(\ZZ)^k$-manifold of which the fixed point set is that of the action $((\ZZ)^k/\ker\rho, C_\rho)$. Then for a fixed point $x$ of the action $((\ZZ)^k, C_\rho)$, we have two tangent $(\ZZ)^k$-representations $(\ZZ)^k\to \GL(T_xM)$ and $(\ZZ)^k\to \GL(T_xC_\rho)$, which induces a $(\ZZ)^k$-representation
\begin{align*}
	(\ZZ)^k \to \GL(T_xM/T_x C_\rho).
\end{align*}
The restriction $\ker\rho\to \GL(T_xM/T_x C_\rho)$ is exactly the normal $\ker\rho$-representation.

\begin{remark} \label{ker_rho_repre}
	Here, we give some necessary explanations and conventions.
	\begin{enumerate}
		\item For any two points $x$ and $y$ in $C_\rho$, the normal $\ker\rho$-representations at $x$ and $y$ are isomorphic, which implies that the tangent $\ker\rho$-representations at $x$ and $y$ in $M$ are identical.
		\item If $C'_\rho$ is another such component of $M^{\ker\rho}$ with the same dimension as $C_\rho$, the normal $\ker\rho$-representations of them may coincide. If so, choose two points $u$ and $v$ from the free parts of $C_\rho$ and $C'_\rho$ under the action of $(\ZZ)^k/\ker\rho$, respectively, and $g_0\in (\ZZ)^k\setminus\ker\rho$. Since the normal $\ker\rho$-representations at $u$ and $v$ are isomorphic, there are open neighborhoods $U$ and $V$ of $M$ at $u$ and $v$, respectively, and a diffeomorphism $F: U\to V$ such that $U$ and $V$ are $\ker\rho$-invariant, and $F$ is $\ker\rho$-equivariant. Then $F$ induces a $\ker\rho$-equivariant diffeomorphism $F': g_0U\to g_0V$. So we can identify the boundaries of $U$ and $V$ via $F$ and identify the boundaries of $g_0U$ and $g_0V$ via $F'$ to receive a new $(\ZZ)^k$-manifold. The condition that $\dim C_\rho = \dim C'_\rho>1$ guarantees that after the equivariant connected sum, $C_\rho$ and $C'_\rho$ become connected. It is clear that the new $(\ZZ)^k$-manifold is equivariantly bordant to $((\ZZ)^k, M)$ since all tangent representations at fixed points are unchanged.
		
		Because of this process, we assume that different components of $M^{\ker\rho}$, with the same dimension larger than one, have distinct normal $\ker\rho$-representations throughout the following. But this is not suitable for two components both of dimension one. Fortunately, this situation is insignificant in our theory.
	\end{enumerate}
\end{remark}

Now return to the construction process. Recall that there is an action of $(\ZZ)^k/\ker\rho\cong \ZZ$ on $C_\rho$ fixing finite points. By a well-known result of \cite{ConnerFloyd}, we know that the fixed point set of this action has even elements.
Therefore, we can always construct a regular connected labeled graph, denoted by $\Gamma_{\rho, C_\rho}$, whose valence is $\dim C_\rho$, whose vertex set is the fixed point set of $((\ZZ)^k/\ker\rho, C_\rho)$, and whose edges are all labeled by $\rho$. In general, the choice of such a graph is not unique. 
If $\dim C_\rho = 1$, $C_\rho$ is equivariantly diffeomorphic to the circle $S^1$ with an involution fixing only two points, and thus 
$\Gamma_{\rho, C_\rho}$ is the graph consisting of two points and one edge joining them. If $\dim C_\rho>1$, there might be many choices for $\Gamma_{\rho, C_\rho}$.

Finally, when $\rho$ and $C_\rho$ run for all possibilities, we have a family of regular connected graphs $\Gamma_{\rho, C_\rho}$. Our desired graph $\Gamma_{((\ZZ)^k, M)}$ is the union of all these graphs $\Gamma_{\rho, C_\rho}$. By Proposition \ref{tangent_representation}, $\Gamma_{((\ZZ)^k, M)}$ is a graph of valence $\dim M=n$.
It is worthwhile to point out that the graph of an $n$-dimensional $(\ZZ)^n$-manifold is unique (i.e., the case of $n=k$), since the action is effective (see also \cite{Lu2009, LuTan2014}).
Moreover, the labeling on edges induces a map
\begin{align*}
	\alpha: E_{\Gamma_{((\ZZ)^k, M)}}\to \Hom_{\ZZ}((\ZZ)^k, \ZZ) \setminus \{0\},
\end{align*}
satisfying that 
\begin{align*}
	\prod_{e\in E_x} \alpha(e)
\end{align*}
is in fact the tangent representation at $x$, where $E_x$ is the set of all edges adjacent to $x$. So 
\begin{align*}
	\phi_*([(\ZZ)^k, M]) = \sum_{x\in M^{(\ZZ)^k}} \prod_{e\in E_x} \alpha(e),
\end{align*}
which implies that the above formula is independent of the choice of $\Gamma_{((\ZZ)^k, M)}$, called the \textit{labeling polynomial of the graph}.

\begin{example} \label{standard_action_real_projective}
	Consider the standard action of $(\ZZ)^n$ on the projective space $\mathbb{R}P^n$ defined by 
	\begin{align*}
		(g_1, \cdots, g_n)^\mathsf{T}\cdot [x_0: x_1: \cdots: x_n] = [x_0: (-1)^{g_1}x_1: \cdots: (-1)^{g_n}x_n].
	\end{align*}
	The fixed points of this action are $x_i = [0: \cdots: 1: \cdots: 0]$ with 1 in the $i$-th place for $0\leqslant i\leqslant n$. A direct computation shows that $\Gamma_{((\ZZ)^n, \mathbb{R}P^n)}$ is the graph in which there are $n+1$ vertices $x_0, \cdots, x_n$ and for any two points $x_i$ and $x_j$ with $i\ne j$, there is an edge joining them, labeling $\rho_i+\rho_j$. That is, $\Gamma_{((\ZZ)^n, \mathbb{R}P^n)}$ is the 1-skeleton of an $n$-simplex.
\end{example}

Hence, we assign a labeled graph to each closed smooth manifold with an effective $(\ZZ)^k$-action fixing a finite set, such that all tangent representations at fixed points are encoded on the labeling of the graph. We give some basic but important properties of this labeled graph.

\begin{proposition}[{\cite{LiLuShen_2025arxiv, Lu2008}}]
	With the preceding notation, $\Gamma_{((\ZZ)^k, M)} = \cup_{\rho, C_\rho} \Gamma_{\rho, C_\rho}$, where each $\Gamma_{\rho, C_\rho}$ is a regular connected subgraph with even vertices and the restriction of $\alpha$ on $E_{\Gamma_{\rho, C_\rho}}$ is constant, satisfying the following three conditions:
	\begin{enumerate}
		\item For each vertex $x$ in $\Gamma_{((\ZZ)^k, M)}$, $\alpha(E_x)$ spans $\Hom_{\ZZ}((\ZZ)^k, \ZZ)$. In particular, $n\geqslant k$.
		\item For each subgraph $\Gamma_{\rho, C_\rho}$, $\alpha(E_x) \equiv \alpha(E_y)\mod\rho$, for all vertices $x$ and $y$ in $\Gamma_{\rho, C_\rho}$.
		\item If $\Gamma_{\rho, C_\rho}$ and $\Gamma_{\rho, C'_\rho}$ are two different subgraphs with the same valence greater than one, then $V_{\Gamma_{\rho, C_\rho}}\cap V_{\Gamma_{\rho, C'_\rho}} = \varnothing$, and $\alpha(E_x)\not\equiv \alpha(E_y)\mod\rho$, for any vertices $x$ in $\Gamma_{\rho, C_\rho}$ and $y$ in $\Gamma_{\rho, C'_\rho}$.
	\end{enumerate}
\end{proposition}

Note that the notation $\alpha(E_x)$ here is not the commonly used symbol for the image set of $E_x$ under $\alpha$. In this paper, $\alpha(E_x)$ is a multiset whose underlying set is $\{\alpha(e): e\in E_x\}$ and for each $e\in E_x$, the multiplicity of $\alpha(e)$ is the number of edges in $E_x$ with the same label as $e$.
The property (1) follows from the effectiveness of the action, and the other two properties have been explained in Remark \ref{ker_rho_repre}.

We emphasize that the labeled graph as a carrier recodes the total information of fixed point data of the action, but what we are more concerned with is the labeling polynomial of the graph and how to detect whether a polynomial in $R_*((\ZZ)^k)$ belongs to the image of $\phi_*$. Li, L\"u and Shen provide a criterion by investigating when the labeling polynomial of a graph satisfying properties (1)-(3) in the above proposition is in the image of $\phi_*$.

\begin{theorem}[{\cite[Theorem 4.1]{LiLuShen_2025arxiv}}] \label{characterization_image}
	Let $k, m$ and $n$ be positive integers with $n\geqslant k$, and $\mathcal{A} = \{\tau_1, \cdots, \tau_m\}$ be a collection of faithful $(\ZZ)^k$-representation classes in $R_n((\ZZ)^k)$. Then the necessary and sufficient condition that $\tau_1+\cdots +\tau_m\in \im\phi_*$ is that for each $\rho\in \Hom_{\ZZ}((\ZZ)^k, \ZZ)$, the set $\{\tau\in \mathcal{A}: \rho \mid \tau\}$ can be decomposed as a disjoint union of nonempty subsets
	\begin{align*}
		\mathcal{A}_\rho^{n_1} \sqcup\cdots \sqcup \mathcal{A}_\rho^{n_l},
	\end{align*}
	with integers $1\leqslant n_1\leqslant \cdots \leqslant n_l$, satisfying the following conditions.
	\begin{enumerate}
		\item For all $1\leqslant i\leqslant l$,
		\begin{enumerate}[label=(1\alph*)]
			\item each representation class of $\mathcal{A}_\rho^{n_i}$ has exactly $n_i$ factors $\rho$;
			\item all representation classes of $\mathcal{A}_\rho^{n_i}$ are isomorphic when restricted to $\ker\rho$;
			\item for any multiset $S$ (including the empty set) containing irreducible $(\ZZ)^k$-representations with cardinality at most $n_i-1$,
			\label{condition_enum_c}
			\begin{align*}
				\sum_{\tau\in \mathcal{A}_\rho^{n_i}} I_\tau(S)\equiv 0\mod 2,
			\end{align*}
			where $I_\tau(S)$ refers to the multiplicity of $S$ in $\tau$ (see \cite[Definition 2.1]{LiLuShen_2025arxiv}).
		\end{enumerate}
		
		\item For all $i\ne j$ and $n_i = n_j>1$, when restricted to $\ker\rho$, each $(\ZZ)^k$-representation of $\mathcal{A}_\rho^{n_i}$ is not isomorphic to any of $\mathcal{A}_\rho^{n_j}$.
	\end{enumerate}
\end{theorem}

\begin{remark}
	If $S$ is the empty set, $I_\tau(S) = 1$, and so the condition \ref{condition_enum_c} implies that each $\mathcal{A}_\rho^{n_i}$ has even elements.
\end{remark}

This necessary and sufficient condition has been proven to be tremendously utile in the argument of the case $n=4$ and $k=3$. On one hand, it is obvious that for each one-dimensional real $(\ZZ)^k$-representation $\rho$, such a decomposition of $\mathcal{A}$ in the theorem is unique. 
This means that there is an algorithm to figure out whether the formal sum $\sum_{\tau\in \mathcal{A}} \tau$ is in the image of $\phi_*$. On the other hand, given a faithful $(\ZZ)^k$-representation class $\tau$ of degree $n$, one can construct a polynomial in $\phi_*(\Z_n((\ZZ)^k))$ whose support contains $\tau$, or show that $\tau$ does not occur in the support of any polynomial in $\phi_*(\Z_n((\ZZ)^k))$, which will be widely used in the next section. 
Here the {\em support} of a polynomial $f$ is the set of all monomials with nonzero coefficients, denoted by $\supp f$.
Thereafter, we keep the notation in this theorem.

\section{Dimension of \texorpdfstring{$\Z_5((\ZZ)^3)$}{}} \label{section_dimension}
Now we come to the main calculation of the work. Before that, some definitions and remarks are needed to be established to simplify the computation.

Let $G$ be a group, $\varphi: G\times M\to M$ be a $G$-action on a manifold $M$, and $\sigma$ be an automomorphism of $G$. Recall that $\varphi_\sigma$ is an action of $G$ on $M$ defined via
\begin{align*}
	\varphi_\sigma(g, x) = \varphi(\sigma(g), x),
\end{align*}
where $g\in G$ and $x\in M$, called the \textit{$\sigma$-induced action of $\varphi$}. A little thought reveals that $(\varphi_\sigma)_{\sigma'} = \varphi_{\sigma\sigma'}$ for automomorphisms $\sigma$ and $\sigma'$ of $G$, and if $\sigma$ is an automorphism of $G$, the effectiveness of $\varphi$ and $\varphi_\sigma$ are equivalent, and they have the same fixed point set.

Returning to the case of $G = (\ZZ)^k$. Let $\varphi: (\ZZ)^k\times M^n\to M^n$ be an effective action with finite fixed points. We explore the relationship of the corresponding polynomials in $R_*((\ZZ)^k)$ between the original action $\varphi$ and the $\sigma$-induced action $\varphi_\sigma$ for all automorphisms $\sigma$ of $(\ZZ)^k$. Write $f = \phi_*[\varphi, M^n]$ and $f_\sigma = \phi_*[\varphi_\sigma, M^n]$. 
If $(\Gamma, \alpha)$ is a $(\ZZ)^k$-labeled graph of $\varphi$, then $(\Gamma, \alpha_\sigma)$ is a $(\ZZ)^k$-labeled graph of $\varphi_\sigma$, where $\alpha_\sigma: E_\Gamma\to \Hom_{\ZZ}((\ZZ)^k, \ZZ)$ is given by
\begin{align*}
	\alpha_\sigma(e)(g) = \alpha(e)(\sigma(g)),
\end{align*}
for all $e\in E_\Gamma$ and $g\in (\ZZ)^k$ (see \cite{lu2010graphsandz2kactions}). It implies that 
\begin{align*}
	f_\sigma = \sum_{x\in M^{(\ZZ)^k}} \prod_{e\in E_x} (\alpha(e)\circ \sigma).
\end{align*}
Moreover, if $x\in M^n$ is a fixed point of the action $\varphi$, since the action is effective, then there always exists an automorphism $\sigma$ of $(\ZZ)^k$ such that the tangent representation at $x$ under the $\sigma$-induced action has the form
\begin{align*}
	\prod_{e\in E_x} \alpha_\sigma(e) = \rho_1\cdots\rho_k \gamma_{k+1}\cdots \gamma_n,
\end{align*}
where $\rho_i$'s are irreducible representations of $(\ZZ)^k$ defined in Section \ref{section_ConnerFloyd_repre} and $\gamma_j$'s are irreducible representations of $(\ZZ)^k$.  

The first step of calculation is to discuss all possibilities of a nonzero monomial in a polynomial of $\phi_*(\Z_n((\ZZ)^k))$. 

\begin{lemma}
	For integers $k\geqslant 2$ and $n\geqslant 2k-1$, the support of every polynomial in $\phi_*(\Z_n((\ZZ)^k))$ does not contain any monomial of the form $\rho_1^{n-k+1}\rho_2\cdots\rho_k$.
\end{lemma}

\begin{proof}
	For an arbitrary $f\in \phi_*(\Z_n((\ZZ)^k))$, if $f=0$, the statement automatically holds. If $f\ne 0$, take $\mathcal{A}=\supp f$. Without loss of generality, assume that $\rho_1^{n-k+1}\rho_2\cdots\rho_k\in \mathcal{A}$. By Theorem \ref{characterization_image}, $\rho_1^{n-k+1}\rho_2\cdots\rho_k$ must lie in a subset of $\mathcal{A}$, say $\mathcal{A}_{\rho_1}^{n-k+1}$. Since $n-k+1\geqslant n-(2k-2)+1\geqslant 2$ and $1\leqslant k-1 \leqslant n-k$, the set $\{\rho_2, \cdots, \rho_k\}$ satisfies the property (1c) in Theorem \ref{characterization_image}, which demonstrates that the number of $\rho_1^{n-k+1}\rho_2\cdots\rho_k$ in $\mathcal{A}_\rho^{n-k+1}$ is even. But it is impossible, and thus the lemma is proved.
\end{proof}

Of course, the above result is true for $k=3$ and $n=5$. Now, let us make a further observation. To simplify notation, set $\rho_{i_1 i_2\cdots i_m}$ to be $\rho_{i_1}+\rho_{i_2}+\cdots \rho_{i_m}\in \Hom_{\ZZ}((\ZZ)^k, \ZZ)$ for integer $m\geqslant 1$ and $1\leqslant i_1< i_2<\cdots< i_m\leqslant k$.

\begin{lemma} \label{not_all_once}
	In the support of every nonzero polynomial in $\phi_*(\Z_5((\ZZ)^3))$, there exists a monomial of the form $\rho_1^2\rho_2\rho_3\gamma$, where $\gamma\in \Hom_{\ZZ}((\ZZ)^3, \ZZ)$.
\end{lemma}

\begin{proof}
	Assume the statement fails for some nonzero polynomial $f\in \phi_*(\Z_5((\ZZ)^3))$, together with the previous lemma, we know that each monomial in $\supp f$ is of the form $\rho_1\rho_2\rho_3\gamma_1\gamma_2$, where $\gamma_1\ne \gamma_2\in \Hom_{\ZZ}((\ZZ)^3, \ZZ)\setminus\{\rho_1, \rho_2, \rho_3\}$. Without loss of generality, suppose that $\rho_1\rho_2\rho_3\rho_{12}\rho_{13}\in \supp f$ or $\rho_1\rho_2\rho_3\rho_{12}\rho_{123}\in \supp f$. Let $\mathcal{A}=\supp f$.
	Then we proceed as follows. 
    
	(1) If $\rho_1\rho_2\rho_3\rho_{12}\rho_{13}\in \supp f$, since $\rho_1\rho_2\rho_3\rho_{12}\rho_{13}$ is isomorphic to $\rho_1\rho_2^2\rho_3^2$ when restricted to $\ker\rho_1$, there is a monomial in $\supp f$ of the form $\gamma_1^2\gamma_2\gamma_3\gamma_4$, where $\gamma_i\in \Hom_{\ZZ}((\ZZ)^3, \ZZ)$.
	
	(2) If $\rho_1\rho_2\rho_3\rho_{12}\rho_{123}\in \supp f$, then $\supp f$ must contain one of the following monomials: $\rho_1\rho_2\rho_3\rho_{12}\rho_{23}$, $\rho_1\rho_2\rho_{13}\rho_{12}\rho_{123}$ or $\rho_1\rho_2\rho_{13}\rho_{12}\rho_{23}$. These three monomials have the same property as $\rho_1\rho_2\rho_3\rho_{12}\rho_{13}$:
	\begin{align*}
		& \rho_1\rho_2\rho_3\rho_{12}\rho_{23} \text{ is isomorphic to } \rho_1^2\rho_2\rho_3^2 \text{ when restricted to } \ker \rho_2; \\
		& \rho_1\rho_2\rho_{13}\rho_{12}\rho_{123} \text{ is isomorphic to } \rho_1^2\rho_2\rho_{13}^2 \text{ when restricted to } \ker \rho_2; \\
		&\rho_1\rho_2\rho_{13}\rho_{12}\rho_{23} \text{ is isomorphic to } \rho_1^2\rho_{13}^2\rho_{12} \text{ when restricted to } \ker \rho_{12}.
	\end{align*}
	These will force $\supp f$ to contain a monomial of the form $\gamma_1^2\gamma_2\gamma_3\gamma_4$.
	
	Hence, both situations contradict the conditions that $f$ satisfies, and the assumption fails.
\end{proof}

The aforementioned discussion can be concluded that in the support of each nontrivial polynomial of $\phi_*(\Z_5((\ZZ)^3))$, there exists a monomial of the form $\rho_1^2\rho_2\rho_3\gamma$, where $\gamma=\rho_2, \rho_{12}, \rho_{23}$ or $\rho_{123}$.

Next, for each $\gamma$, we shall find a polynomial in $\phi_*(\Z_5((\ZZ)^3))$ whose support contains the monomial $\rho_1^2\rho_2\rho_3\gamma$. The search process involves simple but heavy classification discussions, and we just list the results in the following lemma, which can be checked directly by Theorem \ref{characterization_image}.

\begin{lemma} \label{four_poly}
	The four polynomials 
	\begin{align*}
		f_1 &= \rho_1^2\rho_2\rho_3\rho_{123} + \rho_1^2\rho_2\rho_{13}\rho_{23} + \rho_1^2\rho_{12}\rho_3\rho_{23} + \rho_1^2\rho_{12}\rho_{13}\rho_{123}, \\ 
		f_2 &= \rho_1^2\rho_2^2\rho_3 + \rho_1^2\rho_{12}^2\rho_3 + \rho_{12}^2\rho_2^2\rho_3 
		+ \rho_1\rho_{13}\rho_2\rho_{23}\rho_3 + \rho_1\rho_{13}\rho_{12}\rho_{123}\rho_3 + \rho_{12}\rho_{123}\rho_2\rho_{23}\rho_3, \\
		f_3 &= \rho_1^2\rho_2\rho_3\rho_{23} + \rho_1^2\rho_2\rho_{13}\rho_{123} + \rho_1^2\rho_{12}\rho_3\rho_{23} + \rho_1^2\rho_{12}\rho_{13}\rho_{123} + \rho_1\rho_2\rho_3\rho_{12}\rho_{23} \\
		& + \rho_1\rho_2\rho_3\rho_{13}\rho_{23} + \rho_1\rho_2\rho_3\rho_{23}\rho_{123} + \rho_1\rho_2\rho_{12}\rho_{13}\rho_{123} + \rho_1\rho_2\rho_3\rho_{13}\rho_{123} + \rho_1\rho_2\rho_{13}\rho_{23}\rho_{123}, \\
		f_4 &= \rho_1^2\rho_2\rho_{12}\rho_3 + \rho_1^2\rho_2^2\rho_3 + \rho_1^2\rho_2\rho_{12}\rho_{13} + \rho_1^2\rho_{12}^2\rho_{13} + \rho_2^2\rho_1\rho_{12}\rho_3 + \rho_2^2\rho_1\rho_{12}\rho_{23} \\
		& + \rho_2^2\rho_{12}^2\rho_{23} + \rho_{12}^2\rho_1\rho_2\rho_{13} + \rho_{12}^2\rho_1\rho_2\rho_{23} + \rho_1\rho_2\rho_3\rho_{13}\rho_{23} + \rho_1\rho_3\rho_{12}\rho_{13}\rho_{23} + \rho_2\rho_3\rho_{12}\rho_{13}\rho_{23},
	\end{align*}
	are in $\phi_*(\Z_5((\ZZ)^3))$.
\end{lemma}

These four polynomials play an important role in our computation. As we shall see, these four types of polynomials generate the vector space $\phi_*(\Z_5((\ZZ)^3))$ over $\ZZ$. Before proving this, we provide a preliminary lemma.

\begin{lemma} \label{>=4lemma} 
	Let $f\in \phi_*(\Z_n((\ZZ)^k))$ and $\mathcal{A} = \supp f$. If $\tau = \gamma_0^m\gamma_1\cdots\gamma_{n-m}\in \mathcal{A}$, where $m\geqslant 2$ is an integer and $\gamma_i$'s are different irreducible $(\ZZ)^k$-representations, then the subset $\mathcal{A}_{\gamma_0}^m$ of $\mathcal{A}$ containing $\tau$ has more than three elements.
\end{lemma}

\begin{proof}
	By Theorem \ref{characterization_image} (1c), the number of elements in $\mathcal{A}_{\gamma_0}^m$ must be even. If this number is two, then the other element in $\mathcal{A}_{\gamma_0}^m$ is still $\tau$, based on the theorem \ref{characterization_image} (1c) and the fact that $m\geqslant 2$. Therefore, the number $|\mathcal{A}_{\gamma_0}^m|\geqslant 4$.
\end{proof}

Now we are ready for the proof of the main result of this section.

\begin{theorem} \label{generating_set}
	Use the notation in Lemma \ref{four_poly}. The subset 
	\begin{align*}
		\mathcal{P} = \{f_{i, \sigma}: i = 1, 2, 3, 4, \text{ and } \sigma \text{ is an automorphism of } (\ZZ)^3\}
	\end{align*}
	is a generating set of the vector space $\phi_*(\Z_5((\ZZ)^3))$ over $\ZZ$.
\end{theorem}

\begin{proof}
	For every polynomial $f\in \phi_*(\Z_5((\ZZ)^3))$, denote by $s(f)$ the number of monomials in $\supp f$ of the form $\rho_1^2\rho_2\rho_3\gamma$. We will show that every polynomial $f$ in $\phi_*(\Z_5((\ZZ)^3))$ is a linear combination of elements of $\mathcal{P}$ by induction on $s(f)$. 
	
	The base case of induction is satisfied, because a polynomial $f\in \phi_*(\Z_5((\ZZ)^3))$ with $s(f) = 0$ must be trivial by Lemma \ref{not_all_once}.
	
	We assume that any polynomial $f\in \phi_*(\Z_5((\ZZ)^3))$ with $s(f)\leqslant m-1$ is a linear combination of elements of $\mathcal{P}$. Let us consider an arbitrary polynomial $f\in \phi_*(\Z_5((\ZZ)^3))$ with $s(f) = m\geqslant 1$. We have mentioned that there is an automorphism $\sigma$ of $(\ZZ)^3$ such that the support of $f_\sigma$ must include one of the following elements: 
	\begin{align*}
		\rho_1^2\rho_2^2\rho_3, \rho_1^2\rho_2\rho_3\rho_{12}, \rho_1^2\rho_2\rho_3\rho_{23}, \rho_1^2\rho_2\rho_3\rho_{123}.
	\end{align*}
	If $f_\sigma$ is a linear combination of elements of $\mathcal{P}$, so is $f$. Thus without loss of generality, we can assume that one of $\rho_1^2\rho_2^2\rho_3, \rho_1^2\rho_2\rho_3\rho_{12}, \rho_1^2\rho_2\rho_3\rho_{23}, \rho_1^2\rho_2\rho_3\rho_{123}$ is in the support of $f$. The remaining argument will be divided into these four cases. Let $\mathcal{A} = \supp f$.
	
	\begin{enumerate}[label=\textbf{Case (\arabic*).}]
		\setlength{\itemindent}{2.7em}
		\item Suppose that $\rho_1^2\rho_2\rho_3\rho_{123}\in \mathcal{A}$. All five-dimensional faithful representation classes isomorphic to $\rho_1^2\rho_2\rho_3\rho_{123}$ when restricted to $\ker\rho_1$ are 
		\begin{equation} \label{temp_equ1}
			\begin{split}
				\rho_1^2\rho_2\rho_3\rho_{123}, \rho_1^2\rho_2\rho_{13}\rho_{23}, \rho_1^2\rho_{12}\rho_3\rho_{23}, \rho_1^2\rho_{12}\rho_{13}\rho_{123}, \\ \rho_1^2\rho_2\rho_3\rho_{23}, \rho_1^2\rho_2\rho_{13}\rho_{123}, 
				\rho_1^2\rho_{12}\rho_3\rho_{123}, \rho_1^2\rho_{12}\rho_{13}\rho_{23}.
			\end{split}
		\end{equation}
		By Lemma \ref{>=4lemma}, the subset $\mathcal{A}_{\rho_1}^2$ containing $\rho_1^2\rho_2\rho_3\rho_{123}$ has more than three elements. 
		
		\begin{enumerate}[label=\text{Case (1\alph*).}]
			\setlength{\itemindent}{2.7em}
			\item If $\rho_1^2\rho_2\rho_3\rho_{123}$, $\rho_1^2\rho_2\rho_{13}\rho_{23}$, $\rho_1^2\rho_{12}\rho_3\rho_{23}$, $\rho_1^2\rho_{12}\rho_{13}\rho_{123}\in \mathcal{A}_{\rho_1}^2$, take $g = f + f_1 \in \phi_*(\Z_5((\ZZ)^3))$, since $s(g) = s(f+f_1) < s(f) = m$, we know that $g$ is a linear combination of elements of $\mathcal{P}$ by the induction hypothesis, and so is $f = g+f_1$. In a similar way as above, we can always cancel out some elements in $\supp f$ and then decrease $s(f)$ by using elements in $\mathcal{P}$.
			
			\item If $\rho_1^2\rho_2\rho_3\rho_{123}, \tau_1, \tau_2, \tau_3\in \mathcal{A}_{\rho_1}^2$, then we can use polynomials $f_{1, \sigma}$ and $f_{3, \sigma'}$ for certain $\sigma$ and $\sigma'$ to reduce $s(f)$, where $\tau_1$ is an element in the first row of (\ref{temp_equ1}), and $\tau_2$ and $\tau_3$ are elements in the second row of (\ref{temp_equ1}).
		\end{enumerate}
		
		\item The situation that $\rho_1^2\rho_2\rho_3\rho_{23}\in \mathcal{A}$ is similar to the first case.
		
		\item Suppose that $\rho_1^2\rho_2^2\rho_3\in \mathcal{A}$. All five-dimensional representation classes isomorphic to $\rho_1^2\rho_2^2\rho_3$ when restricted to $\ker\rho_1$ are 
		\begin{equation*}
			\rho_1^2\rho_2^2\rho_3, \rho_1^2\rho_2\rho_{12}\rho_3, \rho_1^2\rho_{12}^2\rho_3, \rho_1^2\rho_2^2\rho_{13}, \rho_1^2\rho_2\rho_{12}\rho_{13}, \rho_1^2\rho_{12}^2\rho_{13}.
		\end{equation*}
		If $\rho_1^2\rho_{12}^2\rho_3, \rho_{12}^2\rho_2^2\rho_3\in \mathcal{A}$, let $g = f+f_2\in \phi_*(\Z_5((\ZZ)^3))$, then $s(g) = s(f+f_2) < s(f)$, and so the result is true in this case by the induction hypothesis. If not, then $\rho_1^2\rho_2\rho_{12}\rho_3\in \mathcal{A}$, and the further argument can be reduced to the following case.
		
		\item Suppose that $\rho_1^2\rho_2\rho_3\rho_{12}\in \mathcal{A}$. Also, the subset $\mathcal{A}_{\rho_1}^2$ containing $\rho_1^2\rho_2\rho_3\rho_{12}$ has more than three elements. All five-dimensional representation classes equivalent to $\rho_1^2\rho_2\rho_3\rho_{12}$ when restricted to $\ker\rho_1$ are 
		\begin{equation*}
			\rho_1^2\rho_2\rho_3\rho_{12}, \rho_1^2\rho_2^2\rho_3, \rho_1^2\rho_{12}^2\rho_3, \rho_1^2\rho_2\rho_{13}\rho_{12}, \rho_1^2\rho_2^2\rho_{13}, \rho_1^2\rho_{12}^2\rho_{13}.
		\end{equation*}
		Thus $\mathcal{A}_{\rho_1}^2$ has four elements and $\rho_1^2\rho_2\rho_{13}\rho_{12}\in \mathcal{A}_{\rho_1}^2$. There are two possibilities for the remaining two elements in $\mathcal{A}_{\rho_1}^2$:
		\begin{align*}
			\{\rho_1^2\rho_2^2\rho_3, \rho_1^2\rho_{12}^2\rho_{13}\}\subseteq \mathcal{A}_{\rho_1}^2 \text{ or } \{\rho_1^2\rho_2^2\rho_3, \rho_1^2\rho_2^2\rho_{13}\}\subseteq \mathcal{A}_{\rho_1}^2.
		\end{align*}
		
        When $\{\rho_1^2\rho_2^2\rho_3, \rho_1^2\rho_{12}^2\rho_{13}\}\subseteq \mathcal{A}_{\rho_1}^2$, consider the subset $\mathcal{A}_{\rho_2}^2$ containing $\rho_1^2\rho_2^2\rho_3$ and the subset $\mathcal{A}_{\rho_{12}}^2$ containing $\rho_1^2\rho_{12}^2\rho_{13}$. We perform our arguments as follows.
        If $\rho_{12}^2\rho_2^2\rho_3\in \mathcal{A}_{\rho_2}^2\cup \mathcal{A}_{\rho_{12}}^2$, take $g = f+f_2\in \phi_*(\Z_5((\ZZ)^3))$, since $s(g) = s(f+f_2)<s(f)$, this case is done by the induction hypothesis. 
        Similarly,  if $\rho_2^2\rho_1\rho_{12}\rho_3$, $\rho_2^2\rho_1\rho_{12}\rho_{23}$, $\rho_2^2\rho_{12}^2\rho_{23}$, $\rho_{12}^2\rho_1\rho_2\rho_{13}$, $\rho_{12}^2\rho_1\rho_2\rho_{23}\in \mathcal{A}_{\rho_2}^2\cup \mathcal{A}_{\rho_{12}}^2$, let $g = f+f_4$, we have $s(g) = s(f+f_4)<s(f)$. 
        If $\rho_2^2\rho_1\rho_{12}\rho_3$, $\rho_2^2\rho_1\rho_{12}\rho_{23}$, $\rho_2^2\rho_1^2\rho_{23}$, $\rho_{12}^2\rho_1\rho_2\rho_{13}$, $\rho_{12}^2\rho_1\rho_2\rho_{23}$, $\rho_{12}^2\rho_1^2\rho_{23}\in \mathcal{A}_{\rho_2}^2\cup \mathcal{A}_{\rho_{12}}^2$, let $g = f+f_4$, we have $s(g) = s(f+f_4)<s(f)$. 
        The final situation is that $\rho_2^2\rho_1\rho_{12}\rho_3$, $\rho_2^2\rho_1\rho_{12}\rho_{23}$, $\rho_2^2\rho_1^2\rho_{23}$, $\rho_{12}^2\rho_1\rho_2\rho_{13}$, $\rho_{12}^2\rho_1\rho_2\rho_{23}$, $\rho_{12}^2\rho_2^2\rho_{23}\in \mathcal{A}_{\rho_2}^2\cup \mathcal{A}_{\rho_{12}}^2$, there are five equivalent representation classes when restricted to $\ker\rho_2$, which is impossible.
            
        When $\{\rho_1^2\rho_2^2\rho_3, \rho_1^2\rho_2^2\rho_{13}\}\subseteq \mathcal{A}_{\rho_1}^2$, one can further consider the subset $\mathcal{A}^2_{\rho_2}$ to finish our argument in a similar manner.
	\end{enumerate}
	Together with the above arguments, a polynomial $f\in \phi_*(\Z_5((\ZZ)^3))$ with $s(f) = m$ is a linear combination of elements of $\mathcal{P}$. This completes the induction, so the theorem holds.
\end{proof}

Consequently, it follows that a maximal linearly independent subset of the set $\mathcal{P}$ is a basis of the vector space $\phi_*(\Z_5((\ZZ)^3))$, whose dimension is precisely that of $\Z_5((\ZZ)^3)$ by Theorem \ref{Stong_inj}. Before finding such a subset of $\mathcal{P}$, we prove the following lemma. 

\begin{lemma} \label{four poly-auto}
	Use the notation of Lemma \ref{four_poly}. Let $\sigma$ be an automorphism of $(\ZZ)^3$ which is represented by a non-singular matrix $A$ under the natural basis of $(\ZZ)^3$.
	\begin{enumerate}
		\item $f_{1, \sigma} = f_1$ if and only if $A$ is of the form 
		\begin{align} \label{matrixf1}
			\begin{pmatrix}
				1 & 0 & 0 \\ * & * & * \\ * & * & *
			\end{pmatrix}.
		\end{align}
		\item $f_{2, \sigma} = f_2$ if and only if $A$ is of the form 
		\begin{align} \label{matrixf2}
			\begin{pmatrix}
				* & * & 0 \\ * & * & 0 \\ 0 & 0 & 1
			\end{pmatrix}.
		\end{align}
		\item $f_{3, \sigma} = f_3$ if and only if $A$ is of the form 
		\begin{align*}
			\begin{pmatrix}
				1 & 0 & 0 \\ 0 & 1 & 0 \\ * & * & 1
			\end{pmatrix}.
		\end{align*}
		\item $f_{4, \sigma} = f_4$ if and only if $A$ is one of the matrices
		\begin{align} \label{matrixf4}
			\begin{pmatrix}
				1 & 0 & 0 \\ 0 & 1 & 0 \\ 0 & 0 & 1
			\end{pmatrix}, 
			\begin{pmatrix}
				0 & 1 & 0 \\ 1 & 0 & 0 \\ 0 & 0 & 1
			\end{pmatrix}, 
			\begin{pmatrix}
				1 & 0 & 0 \\ 1 & 1 & 0 \\ 1 & 0 & 1
			\end{pmatrix}, 
			\begin{pmatrix}
				1 & 1 & 0 \\ 1 & 0 & 0 \\ 1 & 0 & 1
			\end{pmatrix}, 
			\begin{pmatrix}
				0 & 1 & 0 \\ 1 & 1 & 0 \\ 0 & 1 & 1
			\end{pmatrix},
			\begin{pmatrix}
				1 & 1 & 0 \\ 0 & 1 & 0 \\ 0 & 1 & 1
			\end{pmatrix}.
		\end{align}
	\end{enumerate}
\end{lemma}

\begin{proof}
	Here we only provide the proofs for (1) and (2), and the proofs of (3) and (4) are similar.
	
	On $f_1$, the ``only if '' part is obvious, since $\rho_1\circ\sigma = \rho_1$. On the other hand, we start with some special properties of $f_1$. Each monomial $\rho_1^2\gamma_1\gamma_2\gamma_3$ in $\supp f_1$ satisfies that $\gamma_i\ne 0$, $\rho_1, \gamma_1, \gamma_2, \gamma_3$ are distinct, and $\rho_1+\gamma_1+\gamma_2+\gamma_3 = 0$.
	So $\rho_1, \gamma_i, \gamma_j$ are linearly independent for all $i\ne j$. 
	Then there are 12 pairs of 1-dimensional $(\ZZ)^3$-representations $\{\gamma_i, \gamma_j\}$ such that $\rho_1, \gamma_i, \gamma_j$ are linearly independent and $\rho_1^2\gamma_i\gamma_j(\rho_1+\gamma_i+\gamma_j)\in \supp f_1$.
    But such 12 pairs $\{\gamma_i, \gamma_j\}$ always exist without the latter condition. That is, for any pair of 1-dimensional $(\ZZ)^3$-representations $\{\gamma, \gamma'\}$ satisfying that $\rho_1, \gamma, \gamma'$ are linearly independent, there exists a monomial $\rho_1^2\gamma\gamma'(\rho_1+\gamma+\gamma')\in \supp f_1$. Suppose that a matrix $A$ is of the form (\ref{matrixf1}) and $\rho_1^2\gamma_1\gamma_2(\rho_1+\gamma_1+\gamma_2)\in \supp f_1$. Because $\rho_1, \gamma_1, \gamma_2$ are linearly independent, $\rho_1\circ\sigma=\rho_1, \gamma_1\circ\sigma, \gamma_2\circ\sigma$ are still linearly independent. As we mentioned above, $\rho_1^2(\gamma_1\circ\sigma)(\gamma_2\circ\sigma)((\rho_1+\gamma_1+\gamma_2)\circ\sigma)\in \supp f_1$. Hence, $f_{1, \sigma} = f_1$.
	
	On $f_2$, if $f_{2, \sigma} = f_2$, it implies that $\rho_3\circ \sigma = \rho_3$ and $\{\rho_1\circ\sigma, \rho_2\circ\sigma, \rho_{12}\circ\sigma\} = \{\rho_1, \rho_2, \rho_{12}\} $. We see that $A$ can only be of the form (\ref{matrixf2}). Conversely, if $A$ is of the form (\ref{matrixf2}), the fact that $f_{2, \sigma} = f_2$ is manifest. 
\end{proof}

With this lemma, it is more convenient to write down all polynomials $f_{i, \sigma}$ for $i=1, 2, 3, 4$ and $\sigma$. For instance, we have already seen that $f_1$ is unchanged under automorphisms of the form (\ref{matrixf1}). Then for an automorphism $\sigma$ represented by the matrix
\begin{align*}
	A = \begin{pmatrix}
		0 & 1 & 0 \\ a_1 & b_1 & c_1 \\ a_2 & b_2 & c_2
	\end{pmatrix},
\end{align*}
since $A = A_1A_2$, where 
\begin{align*}
	A_1 = \begin{pmatrix}
		1 & 0 & 0 \\ b_1 & a_1 & c_1 \\ b_2 & a_2 & c_2
	\end{pmatrix} \text{ and }
	A_2 = \begin{pmatrix}
		0 & 1 & 0 \\ 1 & 0 & 0 \\ 0 & 0 & 1
	\end{pmatrix},
\end{align*}
we know that $f_{1, \sigma} = f_{1, A_1A_2} = (f_{1, A_1})_{A_2} = f_{1, A_2}$. So there are exactly $2^3-1 = 7$ polynomials in the set $\{f_{1, \sigma}: \sigma \text{ is an automorphism of } (\ZZ)^3\}$, each of which can be easily noted down. Taking this discussion as our guide, a maximal linearly independent subset of $\mathcal{P}$ and then the dimension of $\phi_*(Z_5((\ZZ)^3))$ can be determined.

\begin{theorem} \label{thm_dimension}
	The dimension of the vector space $\Z_5((\ZZ)^3)$ over $\ZZ$ is 77.
\end{theorem}

\begin{proof}
	Denote by $\mathcal{P}_i$ the set $\{f_{i, \sigma}: \sigma \text{ is an automorphism of } (\ZZ)^3\}$ for $i=1, 2, 3, 4$. 
	In light of the arguments before this theorem, using Lemma~\ref{four poly-auto}, we first obtain the following observations.
    \begin{enumerate}
        \item $\mathcal{P}_1$ contains $2^3-1 = 7$ elements which are linearly independent. 
        \item $\mathcal{P}_2$ consists of $\frac{(2^3-1)(2^3-2)(2^3-2^2)}{(2^2-1)(2^2-2)} = 28$ elements.
        \item All elements in $\mathcal{P}_1\cup \mathcal{P}_2$ are linearly independent. 
        \item $\mathcal{P}_3$ has $\frac{(2^3-1)(2^3-2)(2^3-2^2)}{4} = 42$ elements. 
        \item $|\mathcal{P}_4| = \frac{(2^3-1)(2^3-2)(2^3-2^2)}{6} = 28$.
    \end{enumerate}
    But not all elements in $\mathcal{P}_1\cup \mathcal{P}_2\cup \mathcal{P}_3$ are linearly independent. 
    By a direct calculation, there are six polynomials with monomials of the form $\rho_1^2\gamma_1\gamma_2\gamma_3$ in $\mathcal{P}_3$: 
	\begin{align*}
		f_{3, 1} &= \rho_1^2\rho_2\rho_3\rho_{23} + \rho_1^2\rho_2\rho_{13}\rho_{123} + \rho_1^2\rho_{12}\rho_3\rho_{23} + \rho_1^2\rho_{12}\rho_{13}\rho_{123} + \rho_1\rho_2\rho_3\rho_{12}\rho_{23} \\
		& + \rho_1\rho_2\rho_3\rho_{13}\rho_{23} + \rho_1\rho_2\rho_3\rho_{23}\rho_{123} + \rho_1\rho_2\rho_{12}\rho_{13}\rho_{123} + \rho_1\rho_2\rho_3\rho_{13}\rho_{123} + \rho_1\rho_2\rho_{13}\rho_{23}\rho_{123}, \\
		f_{3, 2} &= \rho_1^2\rho_3\rho_2\rho_{23} + \rho_1^2\rho_3\rho_{12}\rho_{123} + \rho_1^2\rho_{13}\rho_2\rho_{23} + \rho_1^2\rho_{13}\rho_{12}\rho_{123} + \rho_1\rho_2\rho_3\rho_{13}\rho_{23} \\
		& + \rho_1\rho_2\rho_3\rho_{12}\rho_{23} + \rho_1\rho_2\rho_3\rho_{23}\rho_{123} + \rho_1\rho_3\rho_{12}\rho_{13}\rho_{123} + \rho_1\rho_2\rho_3\rho_{12}\rho_{123} + \rho_1\rho_3\rho_{12}\rho_{23}\rho_{123}, \\
		f_{3, 3} &= \rho_1^2\rho_{12}\rho_3\rho_{123} + \rho_1^2\rho_{12}\rho_{13}\rho_{23} + \rho_1^2\rho_2\rho_3\rho_{123} + \rho_1^2\rho_2\rho_{13}\rho_{23} + \rho_1\rho_{12}\rho_3\rho_2\rho_{123} \\
		& + \rho_1\rho_{12}\rho_3\rho_{13}\rho_{123} + \rho_1\rho_{12}\rho_3\rho_{123}\rho_{23} + \rho_1\rho_{12}\rho_2\rho_{13}\rho_{23} + \rho_1\rho_{12}\rho_3\rho_{13}\rho_{23} + \rho_1\rho_{12}\rho_{13}\rho_{123}\rho_{23}, \\
		f_{3, 4} &= \rho_1^2\rho_{13}\rho_2\rho_{123} + \rho_1^2\rho_{13}\rho_{12}\rho_{23} + \rho_1^2\rho_3\rho_2\rho_{123} + \rho_1^2\rho_3\rho_{12}\rho_{23} + \rho_1\rho_{13}\rho_2\rho_3\rho_{123} \\
		& + \rho_1\rho_{13}\rho_2\rho_{12}\rho_{123} + \rho_1\rho_{13}\rho_2\rho_{123}\rho_{23} + \rho_1\rho_{13}\rho_3\rho_{12}\rho_{23} + \rho_1\rho_{13}\rho_2\rho_{12}\rho_{23} + \rho_1\rho_{13}\rho_{12}\rho_{123}\rho_{23}, \\
		f_{3, 5} &= \rho_1^2\rho_{23}\rho_3\rho_2 + \rho_1^2\rho_{23}\rho_{13}\rho_{12} + \rho_1^2\rho_{123}\rho_3\rho_2 + \rho_1^2\rho_{123}\rho_{13}\rho_{12} + \rho_1\rho_{23}\rho_3\rho_{123}\rho_2 \\
		& + \rho_1\rho_{23}\rho_3\rho_{13}\rho_2 + \rho_1\rho_{23}\rho_3\rho_2\rho_{12} + \rho_1\rho_{23}\rho_{123}\rho_{13}\rho_{12} + \rho_1\rho_{23}\rho_3\rho_{13}\rho_{12} + \rho_1\rho_{23}\rho_{13}\rho_2\rho_{12}, \\
		f_{3, 6} &= \rho_1^2\rho_{123}\rho_3\rho_{12} + \rho_1^2\rho_{123}\rho_{13}\rho_2 + \rho_1^2\rho_{23}\rho_3\rho_{12} + \rho_1^2\rho_{23}\rho_{13}\rho_2 + \rho_1\rho_{123}\rho_3\rho_{23}\rho_{12} \\
		& + \rho_1\rho_{123}\rho_3\rho_{13}\rho_{12} + \rho_1\rho_{123}\rho_3\rho_{12}\rho_2 + \rho_1\rho_{123}\rho_{23}\rho_{13}\rho_2 + \rho_1\rho_{123}\rho_3\rho_{13}\rho_2 + \rho_1\rho_{123}\rho_{13}\rho_{12}\rho_2,
	\end{align*}
    such that $f_{3, 4} = f_{3, 1}+f_{3, 2}+f_{3, 3}$, $f_{3, 5} = f_{3, 2}+f_{3, 3}$, $f_{3, 6} = f_{3, 1}+f_{3, 2}$, and there is no other linearly independence relationship in $\mathcal{P}_1\cup \mathcal{P}_2 \cup \{f_{3, 1}, f_{3, 2}, f_{3, 3}\}$. In other words, up to all automorphisms represented by non-singular matrices of the form 
	\begin{align*}
		\begin{pmatrix}
			1 & 0 & 0 \\ * & * & * \\ * & * & * 
		\end{pmatrix},
	\end{align*}
	there are only three polynomials which, with those polynomials in $\mathcal{P}_1\cup \mathcal{P}_2$, are linearly independent. Likewise, this argument still holds up to automorphisms represented by invertible matrices of the form
	\begin{align*}
		\begin{pmatrix}
			a & b & c \\ * & * & * \\ * & * & * 
		\end{pmatrix},
	\end{align*}
	where $(a, b, c)$ is a nonzero row vector in $(\ZZ)^3$.
	Therefore, a maximal linearly independent subset $S$ of $\mathcal{P}_1\cup \mathcal{P}_2\cup \mathcal{P}_3$ has $7+28+3\times (2^3-1) = 56$ elements.
    
    By a direct calculation, we have that up to all automorphisms of $(\ZZ)^3$ represented by non-singular matrices of the form
	\begin{align*}
		A\begin{pmatrix}
			1 & 0 & 0 \\ 0 & 1 & 0 \\ * & * & 1
		\end{pmatrix},
	\end{align*}
	where $A$ is one of the matrices in (\ref{matrixf4}), there are four polynomials in $\mathcal{P}_4$: 
	\begin{align*}
		f_{4, 1} &= \rho_1^2\rho_2\rho_{12}\rho_3 + \rho_1^2\rho_2^2\rho_3 + \rho_1^2\rho_2\rho_{12}\rho_{13} + \rho_1^2\rho_{12}^2\rho_{13} + \rho_2^2\rho_1\rho_{12}\rho_3 \\ &+ \rho_2^2\rho_1\rho_{12}\rho_{23}
		+ \rho_2^2\rho_{12}^2\rho_{23} + \rho_{12}^2\rho_1\rho_2\rho_{13} + \rho_{12}^2\rho_1\rho_2\rho_{23} + \rho_1\rho_2\rho_3\rho_{13}\rho_{23} \\ &+ \rho_1\rho_3\rho_{12}\rho_{13}\rho_{23} + \rho_2\rho_3\rho_{12}\rho_{13}\rho_{23}, \\
		f_{4, 2} &= \rho_1^2\rho_2\rho_{12}\rho_{13} + \rho_1^2\rho_2^2\rho_{13} + \rho_1^2\rho_2\rho_{12}\rho_3 + \rho_1^2\rho_{12}^2\rho_3 + \rho_2^2\rho_1\rho_{12}\rho_{13} \\ &+ \rho_2^2\rho_1\rho_{12}\rho_{123} 
		+ \rho_2^2\rho_{12}^2\rho_{123} + \rho_{12}^2\rho_1\rho_2\rho_3 + \rho_{12}^2\rho_1\rho_2\rho_{123} + \rho_1\rho_2\rho_{13}\rho_3\rho_{123} \\ &+ \rho_1\rho_{13}\rho_{12}\rho_3\rho_{123} + \rho_2\rho_{13}\rho_{12}\rho_3\rho_{123}, \\
		f_{4, 3} &= \rho_1^2\rho_2\rho_{12}\rho_{23} + \rho_1^2\rho_2^2\rho_{23} + \rho_1^2\rho_2\rho_{12}\rho_{123} + \rho_1^2\rho_{12}^2\rho_{123} + \rho_2^2\rho_1\rho_{12}\rho_{23} \\ &+ \rho_2^2\rho_1\rho_{12}\rho_3
		+ \rho_2^2\rho_{12}^2\rho_3 + \rho_{12}^2\rho_1\rho_2\rho_{123} + \rho_{12}^2\rho_1\rho_2\rho_3 + \rho_1\rho_2\rho_{23}\rho_{123}\rho_3 \\ &+ \rho_1\rho_{23}\rho_{12}\rho_{123}\rho_3 + \rho_2\rho_{23}\rho_{12}\rho_{123}\rho_3, \\
		f_{4, 4} &= \rho_1^2\rho_2\rho_{12}\rho_{123} + \rho_1^2\rho_2^2\rho_{123} + \rho_1^2\rho_2\rho_{12}\rho_{23} + \rho_1^2\rho_{12}^2\rho_{23} + \rho_2^2\rho_1\rho_{12}\rho_{123} \\ &+ \rho_2^2\rho_1\rho_{12}\rho_{13}
		+ \rho_2^2\rho_{12}^2\rho_{13} + \rho_{12}^2\rho_1\rho_2\rho_{23} + \rho_{12}^2\rho_1\rho_2\rho_{13} + \rho_1\rho_2\rho_{123}\rho_{23}\rho_{13} \\ &+ \rho_1\rho_{123}\rho_{12}\rho_{23}\rho_{13} + \rho_2\rho_{123}\rho_{12}\rho_{23}\rho_{13}.
	\end{align*}
    such that $f_{4, 4} = f_{4, 1} + f_{4, 2} + f_{4, 3} + f_{2, 1} + f_{2, 2} + f_{2, 3} + f_{2, 4}$, where 
	\begin{align*}
		f_{2, 1} &= \rho_1^2\rho_2^2\rho_3 + \rho_1^2\rho_{12}^2\rho_3 + \rho_{12}^2\rho_2^2\rho_3 + \rho_1\rho_{13}\rho_2\rho_{23}\rho_3 + \rho_1\rho_{13}\rho_{12}\rho_{123}\rho_3 + \rho_{12}\rho_{123}\rho_2\rho_{23}\rho_3, \\
		f_{2, 2} &= \rho_1^2\rho_2^2\rho_{13} + \rho_1^2\rho_{12}^2\rho_{13} + \rho_{12}^2\rho_2^2\rho_{13} + \rho_1\rho_3\rho_2\rho_{123}\rho_{13} + \rho_1\rho_3\rho_{12}\rho_{23}\rho_{13} + \rho_{12}\rho_{23}\rho_2\rho_{123}\rho_{13}, \\
		f_{2, 3} &= \rho_1^2\rho_2^2\rho_{23} + \rho_1^2\rho_{12}^2\rho_{23} + \rho_{12}^2\rho_2^2\rho_{23} + \rho_1\rho_{123}\rho_2\rho_3\rho_{23} + \rho_1\rho_{123}\rho_{12}\rho_{13}\rho_{23} + \rho_{12}\rho_{13}\rho_2\rho_3\rho_{23}, \\
		f_{2, 4} &= \rho_1^2\rho_2^2\rho_{123} + \rho_1^2\rho_{12}^2\rho_{123} + \rho_{12}^2\rho_2^2\rho_{123} + \rho_1\rho_{23}\rho_2\rho_{13}\rho_{123} 
        \\ &+ \rho_1\rho_{23}\rho_{12}\rho_3\rho_{123} + \rho_{12}\rho_3\rho_2\rho_{13}\rho_{123}
	\end{align*}
	are polynomials in $\mathcal{P}_2$. An easy argument shows that all polynomials in $S\cup \{f_{4, 1}, f_{4, 2}, f_{4, 3}\}$ are linearly independent. 
    This is also true up to automorphisms of $(\ZZ)^3$ represented by invertible matrices of the form 
	\begin{align*}
		A\begin{pmatrix}
			a_1 & b_1 & c_1 \\ a_2 & b_2 & c_2 \\ * & * & *
		\end{pmatrix},
	\end{align*}
	where $A$ is one of the matrices in (\ref{matrixf4}).
    Hence, a maximal linearly independent subset of $\mathcal{P} = \mathcal{P}_1\cup \mathcal{P}_2\cup \mathcal{P}_3 \cup \mathcal{P}_4$ has $56+3\times 7 = 77$ elements, and $\dim_{\ZZ} \Z_5((\ZZ)^3) = \dim_{\ZZ} \phi_*(\Z_5((\ZZ)^3)) = 77$, as desired. 
\end{proof}

So far, the vector space structure of $\phi_*(\Z_5((\ZZ)^3))$ has been fully settled. In the next section, we shall study the preimage of the basis of $\phi_*(\Z_5((\ZZ)^3))$, and construct five-dimensional $(\ZZ)^3$-manifolds whose formal sum of all tangent representations at fixed points are the polynomials in the basis. In particular, once there is a five-dimensional manifold $M$ with an action $\varphi$ of $(\ZZ)^3$ such that $\phi_*[\varphi, M] = f$, then $\phi_*^{-1}(f_\sigma) = [\varphi_\sigma, M]$ for all automorphisms $\sigma$ of $(\ZZ)^3$.
Thus for $i=1, 2, 3, 4$, it suffices to construct an action $\varphi: (\ZZ)^3\times M^5\to M^5$ such that $\phi_*[\varphi, M] = f_{i, \sigma}$ for some $\sigma$.

\section{Representatives of equivariant bordism classes of \texorpdfstring{$\Z_5((\ZZ)^3)$}{}} \label{section_representative}
\subsection{Construction by Milnor hypersurfaces}
First recall the ring structure of the non-equivariant unoriented bordism ring $\Omega_*^O$.

\begin{theorem}[\cite{Milnor1965, Thom1954}] \label{structure_cobordism_ring}
	The unoriented bordism ring $\Omega_*^O$ is the polynomial ring over $\ZZ$ with one variable in every dimension not of the form $2^s-1$. That is, 
	\begin{align*}
		\Omega_*^O\cong \ZZ[x_i: i\ne 2^s-1, s\in \mathbb{N}].
	\end{align*}
	Moreover, if $i$ is even, $x_i$ can be chosen to be the class of the real projective space $\mathbb{R}P^i$; if $i$ is odd and not of the form $2^s-1$, write $i = 2^p(2q+1)-1$ with $p, q\geqslant 1$, then $x_i$ can be selected as the class of the real Milnor hypersurface $H_{2^p, 2^{p+1}q}$, where the real Milnor hypersurface $H_{m, n}$ with $0\leqslant m\leqslant n$ is a submanifold of $\mathbb{R}P^m\times \mathbb{R}P^n$, given via 
	$\{\left([x_0: \cdots: x_m], [y_0: \cdots: y_n]\right)\in \mathbb{R}P^m\times \mathbb{R}P^n: x_0y_0+\cdots +x_my_m=0\}$.
\end{theorem}

Furthermore, from the equivariant case $\Omega_*^{O, (\ZZ)^k}$ to the non-equivariant one, there is a natural homomorphism forgetting the group action
\begin{align*}
	\Omega_*^{O, (\ZZ)^k}\to \Omega_*^O.
\end{align*}
Then a $(\ZZ)^k$-manifold, which does not bound non-equivariantly, does not equivariantly bound (also see \cite{tomDieck1971,KS}). If this action is effective and only has isolated fixed points, it represents a nontrivial element in $\Z_*((\ZZ)^k)$.
In particular, the generators of $\Omega_*^O$ together with appropriate actions of $(\ZZ)^k$ would have the potential to be generators of $\Omega_*^{O, (\ZZ)^k}$, even of $\Z_*((\ZZ)^k)$. 

Tom Dieck \cite{tomDieck1970Fixpunkte} defined actions of $(\ZZ)^k$ on the Milnor hypersurfaces $H_{m, n}$, and Basu, Mukherjee and Sarkar \cite{BasuMukSar2014} considered a special action of $(\ZZ)^n$ on $H_{m, n}$ defined by
\begin{align*}
	&(g_1, \cdots, g_n)^\mathsf{T}\cdot\left([x_0: x_1: \cdots: x_m], [y_0: y_1: \cdots: y_n]\right) \\
	=& \left([x_0: (-1)^{g_1}x_1: \cdots: (-1)^{g_m}x_m], [y_0: (-1)^{g_1}y_1: \cdots: (-1)^{g_m}y_m:\cdots: (-1)^{g_n}y_n]\right),
\end{align*}
which is effective and has only isolated fixed points, and they computed their image under $\phi_*$.
For a positive integer $r$ and subsets $S_1, \cdots, S_n$ of $\{1, 2, \cdots, r\}$, one can define a group homomorphism by
\begin{align*}
	(\ZZ)^r\to (\ZZ)^n, \; T_i \mapsto \sum_{i\in S_j} T_j,
\end{align*}
where $T_i = (0, \cdots, 1, \cdots, 0)^\mathsf{T}$ with $1$ in the $i$-th place. Then the pullback of the $(\ZZ)^n$-action on $H_{m, n}$ via the above group homomorphism induces an action of $(\ZZ)^r$, denoted by $\varphi_{S_1, \cdots, S_n}$.

\begin{proposition}[\cite{BasuMukSar2014}]
	Keep the notation as above.
	\begin{enumerate}
		\item If $S_i$'s are distinct nonempty subsets of $\{1, 2, \cdots, r\}$, then the actions of $(\ZZ)^r$ and $(\ZZ)^n$ on $H_{m, n}$ have the same fixed points.
		\item The image of $[\varphi_{S_1, \cdots, S_n}, H_{m, n}]$ under the homomorphism $\phi_*$ is
		\begin{align*}
			\prod_{i=1}^{m} \rho_{S_i} \sum_{j=1}^n \prod_{k=1, k\ne j}^n \rho_{S_k\vartriangle S_j} + \sum_{i=1}^m \rho_{S_i} \prod_{k=1, k\ne i}^m \rho_{S_k\vartriangle S_i}\left(\prod_{l=1, l\ne i}^n \rho_{S_l} +\sum_{j=1, j\ne i}^n \rho_{S_j}\prod_{l=1, l\ne i, j}^n \rho_{S_l\vartriangle S_j}\right),
		\end{align*}
		where $S_i\vartriangle S_j$ is the symmetric difference of the sets $S_i$ and $S_j$, and $\rho_S$ is the irreducible $(\ZZ)^r$-representation $(g_1, \cdots, g_r)^\mathsf{T}$ such that $g_i = 1$ if $i\in S$ and $0$ otherwise for each subset $S$ of $\{1, \cdots, r\}$.
	\end{enumerate}
\end{proposition}

With their result, we can construct some actions of $(\ZZ)^3$ on the five-dimensional Milnor hypersurfaces with desired tangent representations.

\begin{construction} \label{construct_Milnor_hypersurface}
	Since $5 = 2^1\times (2\times 1+1)-1$, $x_5$ in Theorem \ref{structure_cobordism_ring} can be the class of $H_{2, 4}$. 
	If $S_1, S_2, S_3, S_4$ are distinct nonempty subsets of $\{1, 2, 3\}$, then 
	\begin{equation} \label{image_Milnor_hypersurface}
		\begin{split}
			\phi_*([(\ZZ)^3, H_{2, 4}]) & = \rho_{S_1}\rho_{S_2}\rho_{S_1\vartriangle S_3}\rho_{S_2\vartriangle S_3}\rho_{S_3\triangle S_4} + \rho_{S_1}\rho_{S_2}\rho_{S_1\vartriangle S_4}\rho_{S_2\vartriangle S_4}\rho_{S_3\triangle S_4}\\
			&+ \rho_{S_1}\rho_{S_3}\rho_{S_1\vartriangle S_2}\rho_{S_2\vartriangle S_3}\rho_{S_3\triangle S_4} + \rho_{S_1}\rho_{S_4}\rho_{S_1\vartriangle S_2}\rho_{S_2\vartriangle S_4}\rho_{S_3\triangle S_4}\\
			&+ \rho_{S_2}\rho_{S_3}\rho_{S_1\vartriangle S_2}\rho_{S_1\vartriangle S_3}\rho_{S_3\triangle S_4} + \rho_{S_2}\rho_{S_4}\rho_{S_1\vartriangle S_2}\rho_{S_1\vartriangle S_4}\rho_{S_3\triangle S_4}.
		\end{split}
	\end{equation}
	
	(1) Choose $S_1 = \{2\}$, $S_2 = \{1, 2\}$, $S_3 = \{2, 3\}$ and $S_4 = \{1, 2, 3\}$. It derives an effective action $\varphi_1 = \varphi_{S_1, \cdots, S_4}$ of $(\ZZ)^3$ on $H_{2, 4}$ with finite fixed points such that $\phi_*([\varphi_1, H_{2, 4}]) = f_1$.
	
	(2) Choose $S_1 = \{1\}$, $S_2 = \{2\}$, $S_3 = \{1, 2\}$ and $S_4 = \{1, 2, 3\}$. Then we acquire an effective action $\varphi_2 = \varphi_{S_1, \cdots, S_4}$ of $(\ZZ)^3$ on $H_{2, 4}$ with finite fixed points such that $\phi_*([\varphi_2, H_{2, 4}]) = f_2$.
\end{construction}

By the equation (\ref{image_Milnor_hypersurface}), we cannot find subsets $S_1, S_2, S_3$ and $S_4$ of $\{1, 2, 3\}$ such that 
\begin{align*}
	\phi_*([\varphi_{S_1, S_2, S_3, S_4}, H_{2, 4}]) = f_3 \text{ or } f_4,
\end{align*}
and to remedy this we turn to the next approach.

\subsection{Construction by small covers}
This approach to discover representatives of the equivariant unoriented bordism classes of $\Z_5((\ZZ)^3)$ is inspired by two facts. One is the following lemma.

\begin{lemma}
	For each $\beta\in \Z_n((\ZZ)^k)$ with $0<k\leqslant n$, there exists an action $\varphi$ of $(\ZZ)^n$ on a manifold $M^n$ and a subgroup $H$ of $(\ZZ)^n$ of rank $k$ such that $\beta = [(\ZZ)^k, M^n]$ in $\Z_n((\ZZ)^k)$, where the action of $(\ZZ)^k$ on $M^n$ is the pullback of the action $\varphi$ through the group homomorphism $(\ZZ)^k\cong H \hookrightarrow(\ZZ)^n$.
\end{lemma}

\begin{proof}
	Let $((\ZZ)^k, M^n)$ be a representative of $\beta$ which is an effective action with finite fixed points. Identify $(\ZZ)^n$ with $(\ZZ)^k\oplus (\ZZ)^{n-k}$. Define an action of $(\ZZ)^n$ on $M$ by $(g_1, g_2)\cdot x = g_1x$. Let $H = (\ZZ)^k\oplus \{0\}< (\ZZ)^n$. Then clearly $H\cong (\ZZ)^k$ and the subgroup $H$-action on $M$ is exactly the original action of $(\ZZ)^k$ on $M$.
\end{proof}

Therefore, with the above lemma, to find the five-dimensional $(\ZZ)^3$-manifolds with the desired image under the homomorphism $\phi_*$, our strategy is as follows: we first begin with an effective five-dimensional $(\ZZ)^5$-manifold with isolated fixed points, and then determine a subgroup of $(\ZZ)^5$ of rank 3 such that the fixed points of the subgroup action are still isolated. Finally we calculate the image of the action under $\phi_*$.
Consequently, it is essential to figure out how to identify these five-dimensional $(\ZZ)^5$-manifolds.
It is motivated by L\"u and Tan, who gave an affirmative answer in \cite{LuTan2014} to the conjecture: each class of $\Z_n((\ZZ)^n)$ contains a small cover as its representative. 

Small covers are introduced by Davis and Januszkiewicz \cite{smallcover}, indicating a profound connection between equivariant topology and combinatorics of simple convex polytopes. 
An $n$-dimensional \textit{small cover} over a simple convex polytope $P^n$ is an $n$-dimensional manifold $M^n$ with a locally standard action of $(\ZZ)^n$ whose orbit space is homeomorphic to $P^n$, often denoted by $\pi: M^n\to P^n$, where $\pi$ is the orbit map. Here, a locally standard action of $(\ZZ)^n$ on an $n$-dimensional manifold means that the action is locally isomorphic to the standard action of $(\ZZ)^n$ on $\mathbb{R}^n$. The real projective space $\mathbb{R}P^n$ with action in Example \ref{standard_action_real_projective} is a small cover over an $n$-simplex.
Clearly, the action in a small cover is effective and has finite fixed points, and thus the equivariant bordism class of an $n$-dimensional small cover belongs to $\Z_n((\ZZ)^n)$.
Each small cover over $P^n$ corresponds to a function $\lambda$, called the \textit{characteristic function}, assigning a nonzero element in $(\ZZ)^n$ to every codimension-one face of $P^n$, satisfying the following condition.

\begin{condition} \label{characteristic_function}
	If $F_1, \cdots, F_n$ are codimension-one faces of $P^n$ meeting at some vertex, then $\lambda(F_1), \cdots, \lambda(F_n)$ form a basis of $(\ZZ)^n$. 
\end{condition}

Conversely, we know from \cite{smallcover} that for a function $\lambda$ satisfying (\ref{characteristic_function}), there is a unique small cover with $\lambda$ as its characteristic function, up to equivariant homeomorphism.
It indicates that all information on equivariant topology of a small cover is embodied in the corresponding polytope $P$ and its characteristic function $\lambda$. 
For example, the mod 2 cohomology ring of a small cover can be completely described by the pair $(P, \lambda)$ (see \cite[Theorem 4.14]{smallcover}). 
Meanwhile, the characteristic function of a small cover $\pi: M^n\to P^n$ has also a direct connection with $\phi_*([(\ZZ)^n, M^n])$, as we summarize as the following lemma.

\begin{lemma}[\cite{Lu2008, LuTan2014}] \label{characfun_tangentrepre}
	Let $\pi: M^n\to P^n$ be a small cover with the $(\ZZ)^n$-labeled graph $(\Gamma, \alpha)$ and the characteristic function $\lambda$. Then $\Gamma$ is the one-skeleton of $P^n$, and $\alpha$ and $\lambda$ determine each other. Precisely, if $v$ is a vertex of $P^n$, there are $n$ codimension-one faces $F_1, \cdots, F_n$ of $P^n$ and $n$ dimension-one faces $e_1, \cdots, e_n$ of $P^n$ such that $v = F_1\cap\cdots\cap F_n = e_1\cap\cdots\cap e_n$. Without loss of generality, assume that $e_i = F_1\cap \cdots\cap F_{i-1}\cap F_{i+1}\cap \cdots \cap F_n$. Then 
	\begin{align*}
		\alpha(e_i) (\lambda(F_j)) = \delta_{ij}.
	\end{align*}
\end{lemma}

Once the five-dimensional $(\ZZ)^5$-manifold $((\ZZ)^5, M^5)$ is obtained, we advance to find a subgroup $H$ of $(\ZZ)^5$ of rank 3 such that the subgroup action only has isolated fixed points. From Corollary \ref{subgroup_tangent_repre}, we know that if $x$ is a fixed point of $((\ZZ)^5, M^5)$, the tangent representation at $x$ in $(H, M^5)$ is clear. In general, if $x\in M^H \setminus M^{(\ZZ)^5}$, the tangent representation at $x$ in $(H, M^5)$ is hard to detect. Hence, it is reasonable to require that $M^H = M^{(\ZZ)^5}$ and to determine when $H$ satisfies this property. 

\begin{lemma}
    Suppose that $G$ is a group acting on a set $X$ and $H$ is a subgroup of $G$. Then $X^H = X^G$ if and only if $H$ does not contain in the isotropy subgroup $G_x$ for all $x\in X \setminus X^G$.
\end{lemma}

\begin{proof}
    Since $X^G\subseteq X^H$, it suffices to show that $X^H\setminus X^G\ne\varnothing\iff H\subseteq G_x$ for some $x\in X\setminus X^G$, which is clear from the fact that $x\in X^H\iff H\subseteq G_x$.
\end{proof}

\begin{remark}
    When $G$ is a torus $T^k$ and $X$ is an orientable manifold, such a subgroup $H$ always exists. In fact, it was shown in \cite[Theorem IV.10.5]{Bredon_transformationgroup} and \cite[Lemma 7.4.3]{BP} that if $M$ is a connected orientable manifold equipped with an effective action of a torus $T^k$, then there exists a circle subgroup $H\subseteq T^k$ such that $M^{H} = M^{T^k}$.
\end{remark}

Since all isotropy subgroups of a small cover can be determined by the characteristic function, we can answer when a subgroup action fixes the same points for a small cover by use of the characteristic function.

\begin{corollary} \label{subgroup_fix_same}
	Let $\pi: M^n\to P^n$ be a small cover, $\lambda$ be its characteristic function, and $H$ be a subgroup of $(\ZZ)^n$. Then $M^H = M^{(\ZZ)^n}$ if and only if for any $n-1$ codimension-one faces $F_1, \cdots, F_{n-1}$ with $F_1\cap \cdots \cap F_{n-1} \ne \varnothing$, $H\not\subseteq \text{span}\{\lambda(F_1), \cdots, \lambda(F_{n-1})\}$.
\end{corollary}

\begin{proof}
	If $F^r$ is an $r$-face of $P^n$, $1\leqslant r\leqslant n$, and $x\in \pi^{-1}(\text{int} F^r)$, then there are $n-r$ codimension-one faces $F_1, \cdots, F_{n-r}$ of $P^n$ such that $F^r = F_1\cap \cdots\cap F_{n-r}$ and $G_x$ is the subgroup spanned by $\lambda(F_1), \cdots, \lambda(F_{n-r})$ (see \cite{smallcover}). Thus, the result follows from the above lemma.
\end{proof}
 
Based upon the preceding discussion, now we can construct representatives of $\phi_*^{-1}(f_i)$, $i = 3, 4$. Apart from what was mentioned above, L\"u and Tan \cite{LuTan2014} also showed that the vector space $\Z_n((\ZZ)^n)$ is generated by the classes of small covers over $\triangle^{n_1}\times \cdots \times \triangle^{n_l}$, where $n_1+\cdots +n_l = n$ and $\triangle^r$ is an $r$-simplex. So we shall start with small covers over a product of simplexes.

\begin{construction} \label{construct_small_cover1}
	Consider the five-dimensional simple convex polytope $\triangle^1\times \triangle^4$. Let $F^1_i$ be codimension-one faces of $\triangle^1$, and $F^2_j$ be codimension-one faces of $\triangle^4$ for $i=1, 2$ and $j=1, 2, 3, 4, 5$, respectively. Then all codimension-one faces of $\triangle^1\times \triangle^4$ are $F^1_1\times \triangle^4, F^1_2\times \triangle^4, \triangle^1\times F^2_1, \triangle^1\times F^2_2, \triangle^1\times F^2_3, \triangle^1\times F^2_4$ and $\triangle^1\times F^2_5$. Denote the set of these codimension-one faces by $\mathscr{F}$. Let $\lambda: \mathscr{F}\to (\mathbb{Z}_2)^5$ be the function indicated as the following matrix,
	\begin{align*}
		\bordermatrix{
			&F^1_1\times \triangle^4& F^1_2\times \triangle^4& \triangle^1\times F^2_1& \triangle^1\times F^2_2& \triangle^1\times F^2_3& \triangle^1\times F^2_4 &\triangle^1\times F^2_5\cr
			&1&0&1&0&0&0&1\cr
			&0&1&1&0&0&0&1\cr
			&0&0&0&1&0&0&1\cr
			&0&0&0&0&1&0&1\cr
			&0&0&0&0&0&1&1
		},
	\end{align*}
	which satisfies the condition (\ref{characteristic_function}). So there is a locally standard action of $(\ZZ)^5$ on a five-dimensional manifold $M_1$ such that it is a small cover over $\triangle^1\times \triangle^4$ with characteristic function $\lambda$.
	By the lemma \ref{characfun_tangentrepre}, we can compute all the tangent representations of this action, which are
	\begin{align*}
		\rho_2\rho_{23}\rho_{24}\rho_{25}\rho_{12}, 
		\rho_3\rho_{23}\rho_{34}\rho_{35}\rho_{12}, 
		\rho_4\rho_{24}\rho_{34}\rho_{45}\rho_{12},
		\rho_5\rho_{25}\rho_{35}\rho_{45}\rho_{12},
		\rho_2\rho_3\rho_4\rho_5\rho_{12},\\
		\rho_1\rho_{13}\rho_{14}\rho_{15}\rho_{12},
		\rho_3\rho_{13}\rho_{34}\rho_{35}\rho_{12},
		\rho_4\rho_{14}\rho_{34}\rho_{45}\rho_{12},
		\rho_5\rho_{15}\rho_{35}\rho_{45}\rho_{12},
		\rho_1\rho_3\rho_4\rho_5\rho_{12}.
	\end{align*}
	Let $H$ be the subgroup of $(\mathbb{Z}_2)^5$ spanned by
	\begin{align*}
		(0, 1, 1, 1, 1)^\mathsf{T}, (1, 1, 0, 1, 0)^\mathsf{T}, (1, 1, 0, 0, 1)^\mathsf{T}.
	\end{align*}
	One can directly verify that there are no codimension-one faces $F_1, F_2, F_3$ and $F_4$ such that $H$ contains in the group spanned by $\lambda(F_1), \lambda(F_2), \lambda(F_3)$ and $\lambda(F_4)$. Then by Corollary \ref{subgroup_fix_same}, $M_1^H = M_1^{(\mathbb{Z}_2)^5}$. Moreover, 
	\begin{align*}
		H^\perp = \left\{ {\bf 0}, (1, 1, 1, 0, 0)^\mathsf{T}, (1, 0, 0, 1, 1)^\mathsf{T}, (0, 1, 1, 1, 1)^\mathsf{T} \right\}.
	\end{align*}
	So in the restricted action $\varphi_3$ of $H$, the above ten tangent representations become
	\begin{align*}
		\rho_2\rho_1\rho_{125}\rho_{25}\rho_{12}, \rho_{12}^2\rho_1\rho_{25}\rho_{125}, \rho_{15}\rho_{125}\rho_{25}\rho_1\rho_{12}, \rho_5\rho_{25}\rho_{125}\rho_1\rho_{12}, \rho_2\rho_{12}^2\rho_{15}\rho_5,\\ 
		\rho_1\rho_2\rho_5\rho_{15}\rho_{12}, \rho_{12}^2\rho_2\rho_{25}\rho_{125}, \rho_{15}\rho_5\rho_{25}\rho_1\rho_{12}, \rho_5\rho_{15}\rho_{125}\rho_1\rho_{12}, \rho_1\rho_{12}^2\rho_{15}\rho_5,
	\end{align*}
	and thus there is an automorphism $\sigma$ of $(\ZZ)^3$ such that $\phi_*([\varphi_3, M_1]) = f_{3, \sigma}$.
\end{construction}

\begin{construction} \label{construct_small_cover2}
	Now we consider the simple convex polytope $\triangle^2\times \triangle^3$. Let $F^1_i$ be codimension-one faces of $\triangle^2$ and $F^2_j$ be codimension-one faces of $\triangle^4$ for $i = 1, 2, 3$ and $j = 1, 2, 3, 4$, respectively. 
	Denote the set of all codimension-one faces of $\triangle^2\times \triangle^3$ by $\mathscr{F}$. Let $\lambda: \mathscr{F}\to (\mathbb{Z}_2)^5$ be the following function.
	\begin{align*}
		\bordermatrix{
			&F^1_1\times \triangle^3& F^1_2\times \triangle^3 & F^1_3\times \triangle^3 & \triangle^2\times F^2_1& \triangle^2\times F^2_2& \triangle^2\times F^2_3& \triangle^2\times F^2_4 \cr
			&1&0&1&1&0&0&1\cr
			&0&1&1&1&0&0&1\cr
			&0&0&0&1&0&0&1\cr
			&0&0&0&0&1&0&1\cr
			&0&0&0&0&0&1&1
		}
	\end{align*}
	It is also a characteristic function of $M_2\to \triangle^2\times\triangle^3$ for some five-dimensional $(\ZZ)^5$-manifold $M_2$ by direct verification.
	The corresponding tangent representations at $12$ fixed points are 
	\begin{align*}
		&\rho_5\rho_{23}\rho_3\rho_4\rho_{13}, \rho_5\rho_{23}\rho_3\rho_4\rho_{12}, \rho_5\rho_{13}\rho_3\rho_4\rho_{12}, \rho_5\rho_{23}\rho_{35}\rho_{45}\rho_{13}, \\ &\rho_5\rho_{23}\rho_{35}\rho_{45}\rho_{12}, \rho_5\rho_{13}\rho_{35}\rho_{45}\rho_{12}, 
		\rho_4\rho_{23}\rho_{34}\rho_{45}\rho_{13}, \rho_4\rho_{23}\rho_{34}\rho_{45}\rho_{12},\\ &\rho_4\rho_{13}\rho_{34}\rho_{45}\rho_{12}, \rho_3\rho_{23}\rho_{34}\rho_{35}\rho_{13}, \rho_3\rho_{23}\rho_{34}\rho_{35}\rho_{12}, \rho_3\rho_{13}\rho_{34}\rho_{35}\rho_{12}.
	\end{align*}
	Let $H$ be the subgroup of $(\mathbb{Z}_2)^5$ spanned by 
	\begin{align*}
		(0, 1, 1, 1, 1)^\mathsf{T}, (0, 1, 0, 1, 0)^\mathsf{T}, (1, 1, 1, 0, 0)^\mathsf{T}.
	\end{align*}
	There are no codimension-one faces $F_1, F_2, F_3$ and $F_4$ such that $H$ contains in the group spanned by $\lambda(F_1), \lambda(F_2), \lambda(F_3)$ and $\lambda(F_4)$. Then $M_2^H = M_2^{(\mathbb{Z}_2)^5}$. Moreover, 
	\begin{align*}
		H^\perp = \{{\bf 0}, (1, 1, 0, 1, 0)^\mathsf{T}, (1, 0, 1, 0, 1)^\mathsf{T}, (0, 1, 1, 1, 1)^\mathsf{T}\}.
	\end{align*}
	So in the subgroup action $\varphi_4$ of $H$, the above 12 tangent representations change into
	\begin{align*}
		&\rho_{13}^2\rho_{23}\rho_3\rho_{12}, \rho_{13}\rho_{23}\rho_3\rho_{12}^2, \rho_{13}^2\rho_3\rho_{12}^2, \rho_{13}^2\rho_{23}^2\rho_1, \rho_{13}\rho_{23}^2\rho_1\rho_{12}, \rho_{13}^2\rho_1\rho_{23}\rho_{12}, \\
		&\rho_{12}\rho_{23}^2\rho_{123}\rho_{13}, \rho_{12}^2\rho_{23}^2\rho_{123}, \rho_{12}^2\rho_{13}\rho_{123}\rho_{23}, \rho_3\rho_{23}\rho_{123}\rho_1\rho_{13}, \rho_3\rho_{23}\rho_{123}\rho_1\rho_{12}, \rho_3\rho_{13}\rho_{123}\rho_1\rho_{12}.
	\end{align*}
	Then $\phi_*([\varphi_4, M_2]) = f_{4, \sigma}$ for some automorphism $\sigma$ of $(\ZZ)^3$.
\end{construction}

\begin{remark}
	In fact, $M_1$ is the total space of a $\mathbb{R}P^4$-bundle over $\mathbb{R}P^1$, and $M_2$ is the total space of a $\mathbb{R}P^3$-bundle over $\mathbb{R}P^2$ (see Remark \ref{remark_projectivization} for details).
\end{remark}

\begin{remark}
	However, any small cover $M^5\to \triangle^5$ cannot produce the desired representative with the approach outlined earlier. There is essentially only one possible characteristic function $\lambda$, namely the one indicated below,
	\begin{align*}
		\begin{pmatrix}
			1&0&0&0&0&1\\
			0&1&0&0&0&1\\
			0&0&1&0&0&1\\
			0&0&0&1&0&1\\
			0&0&0&0&1&1
		\end{pmatrix}.
	\end{align*}
	So all tangent representations of this action are
	\begin{align*}
		\rho_1\rho_2\rho_3\rho_4\rho_5,
		\rho_1\rho_{12}\rho_{13}\rho_{14}\rho_{15},
		\rho_2\rho_{12}\rho_{23}\rho_{24}\rho_{25},
		\rho_3\rho_{13}\rho_{23}\rho_{34}\rho_{35},
		\rho_4\rho_{14}\rho_{24}\rho_{34}\rho_{45},
		\rho_5\rho_{15}\rho_{25}\rho_{35}\rho_{45}.
	\end{align*}
	If there is a subgroup $H$ of $(\ZZ)^5$ such that $H\cong (\ZZ)^3$, $M^H = M^{(\ZZ)^5}$ and the equivariant bordism class of $(H, M^5)$ is nontrivial, then there must be different integers $i, j$ and $k$ such that $\rho_{ij} + H^\perp= \rho_{ik} + H^\perp$. Then $\rho_{jk} +H^\perp = 0$. It is impossible since $\rho_{jk}$ must appear in the above tangent representation classes and all fixed points are isolated. 
\end{remark}

In conclusion, the group structure of $\Z_5((\ZZ)^3)$ is fully understood by Theorem \ref{Stong_inj}, Theorem \ref{generating_set}, Construction \ref{construct_Milnor_hypersurface}, \ref{construct_small_cover1} and \ref{construct_small_cover2}.

\begin{theorem} \label{thm_representative}
	The vector space $\Z_5((\ZZ)^3)$ over $\ZZ$ is generated by the classes of $(\varphi_{1, \sigma}, H_{2, 4})$, $(\varphi_{2, \sigma}, H_{2, 4})$, $(\varphi_{3, \sigma}, M_1)$, and $(\varphi_{4, \sigma}, M_2)$ for all automorphisms $\sigma$ of $(\ZZ)^3$.
\end{theorem}

\begin{remark} \label{remark_projectivization}
	In \cite{Pergher1993}, Pergher showed that every $(\ZZ)^k$-manifold is equivariantly bordant to the projectivization of some real vector bundle with certain $(\ZZ)^k$-action constructible from the fixed data of the original $(\ZZ)^k$-action. We achieve a similar result in our setting and the theorem \ref{thm_representative} indicates exactly what the projectivization is.
	
	Note that the Milnor hypersurface can be regraded as the total space of a fiber bundle with fiber the projective space \cite{Milnor1965}. Explicitly, suppose that $m\leqslant n$ are nonnegative integers. Denote by $\gamma$ the real line bundle $E\to \mathbb{R}P^m$, where $E = \{(x, y)\in \mathbb{R}P^m\times \mathbb{R}^{n+1}: y\in x\}$. Then $\gamma$ is a subbundle of the trivial bundle $\varepsilon$ of rank $n+1$ over $\mathbb{R}P^m$. Clearly, the Milnor hypersurface $H_{m, n}$ is the projectivization of the orthogonal complement of $\gamma$ in $\varepsilon$.
	
	Also, a small cover over a product of simplices is equivalent in the sense of \cite{smallcover}, in particular, weakly equivariantly homeomorphic, to a generalized real Bott manifold \cite{ChoiMasudaSuh2010}, which is also the projectivization of some vector bundle. Precisely, a generalized real Bott tower is a sequence of projective bundles 
	\begin{align*}
		B_N \to B_{N-1}\to\cdots\to B_1\to B_0 = \{\text{a point}\},
	\end{align*}
	where $B_i$ is the projectivization of the Whitney sum of $n_i+1$ real line bundles over $B_{i-1}$ for $i=1, \cdots, N$, and each $B_i$ is called a generalized real Bott manifold. Every $i$-stage $B_i$ naturally provides a small cover over $\prod_{j=1}^i \triangle^{n_j}$ and vice versa. 
	
	To sum up, every element in $\Z_5((\ZZ)^3)$ has a representative which is the projectivization of a real vector bundle over the disjoint union of the real projective spaces with certain action of $(\ZZ)^3$. 
\end{remark}


\end{document}